\documentclass[OJMO,manuscript]{cedram}

\usepackage{mathtools} 

\usepackage[linesnumbered,boxed,english,onelanguage,vlined,algo2e]{algorithm2e} 

\usepackage{accents}
\newcommand\thickbar[1]{\accentset{\rule{.4em}{.7pt}}{#1}}

\usepackage{multirow}
\DeclareFontFamily{U}{skulls}{}
\DeclareFontShape{U}{skulls}{m}{n}{ <-> skull }{}
\newcommand{\skull}{\text{\usefont{U}{skulls}{m}{n}\symbol{'101}}}

\usepackage{microtype} 
\usepackage{xurl} 

\makeatletter
\renewcommand*{\top}{%
  {\mathpalette\@transpose{}}%
}
\newcommand*{\@transpose}[2]{%
  \raisebox{\depth}{$\m@th#1\scriptscriptstyle\mathsf{T}$}%
}
\makeatother

\newcommand{\newatop}[2]{\genfrac{}{}{0pt}{1}{#1}{#2}}

\DeclareMathOperator{\rank}{rank}
\DeclareMathOperator*{\argmax}{argmax}
\DeclareMathOperator*{\argmin}{argmin}

\DeclareMathOperator{\tr}{\textstyle{tr}}

\title[Row-Sparse
\lowercase{ah}-Symmetric Reflexive Generalized Inverses]{Good and Fast Row-Sparse
\lowercase{ah}-Symmetric\\ Reflexive Generalized Inverses}

\author{\firstname{Gabriel} \lastname{Ponte}}
\address{University of Michigan, Ann Arbor, USA}
\thanks{G. Ponte was supported in part by CNPq GM-GD scholarship 161501/2022-2.}
\email{gabponte@umich.edu}
\author{\firstname{Marcia} \lastname{Fampa}}
\address{Federal University of Rio de Janeiro, Brazil }
\email{fampa@cos.ufrj.br}
\thanks{M. Fampa was supported in part by CNPq grant 307167/2022-4.} 
\author{\firstname{Jon} \lastname{Lee}}
\address{University of Michigan, Ann Arbor, USA}
\thanks{J. Lee was supported in part by AFOSR grant FA9550-22-1-0172.}
\email{jonxlee@umich.edu}
\author{\firstname{Luze} \lastname{Xu}}
\address{University of California, Davis, USA}
\email{lzxu@ucdavis.edu}
\ifpdf
\hypersetup{
  pdftitle={Good and Fast Row-Sparse 
ah-Symmetric Reflexive Generalized Inverses},
  pdfauthor={G. Ponte, M. Fampa, J. Lee, and L. Xu}
}
\fi

  
\begin{abstract} 
We present several algorithms aimed at constructing sparse and structured sparse (row-sparse) generalized inverses, with application to  the efficient computation of least-squares solutions,
for inconsistent systems of linear equations, in the setting of multiple right-hand sides and a rank-deficient constraint matrix. Leveraging our earlier formulations 
to minimize the 1- and 2,1-norms 
of generalized inverses that satisfy important properties of the  Moore-Penrose pseudoinverse, we  develop efficient and scalable ADMM algorithms to address these norm-minimization problems and to limit the number of nonzero rows in the solution. 
We establish a 2,1-norm approximation result for a local-search procedure that was originally designed 
for 1-norm minimization, and we compare the ADMM algorithms  with the local-search procedure and with general-purpose optimization solvers.
\end{abstract}

\begin{document}

\maketitle


\section{Introduction}\label{sec:introd}
The M-P (Moore-Penrose) pseudoinverse is an important object in
numerical linear algebra. It has many uses; in particular, and relevant to our study, it can be used to compute least-squares solutions in the context of linear statistical models.
Given a real singular-value decomposition $A\!=\!U\Sigma V^\top$ 
of $A$, the M-P pseudoinverse of $A$ is $A^{\dagger}\!:=\!V\Sigma^{\dagger} U^\top$, 
where the diagonal matrix $\Sigma^{\dagger}$
has the shape of the transpose of the diagonal matrix $\Sigma$, and is calculated from $\Sigma$
by taking reciprocals of the nonzero elements of $\Sigma$, and is otherwise 0 (see \cite[Section 2.5.3]{GVL1996}, for example).
The following well-known theorem  characterizes
 the M-P pseudoinverse.
\begin{theo}[see \protect{\cite[Theorem 1]{Penrose}}]
For $A\!\in\!\mathbb{R}^{m \times n}$, the M-P pseudoinverse $A^{\dagger}$ is the unique 
 $H\!\in\!\mathbb{R}^{n \times m}$ \hbox{satisfying:}
	\begin{align}
\mbox{(generalized inverse)\quad	 }	& AHA = A \label{property1} \tag{P1}\\
\mbox{(reflexive)\quad }		& HAH = H \label{property2} \tag{P2}\\
\mbox{(ah-symmetric)\quad }		& (AH)^{\top} = AH \label{property3} \tag{P3}\\
\mbox{(ha-symmetric)\quad }		& (HA)^{\top} = HA \label{property4} \tag{P4}
	\end{align}
\end{theo}

As is common (see \cite[Introduction, Section 2]{BenIsrael1974}),
we say that a \emph{generalized inverse} of $A$ is any
$H$ satisfying \ref{property1}, and a generalized inverse of $A$ is \emph{reflexive} if it additionally satisfies \ref{property2}.
Following \cite[p. 1723]{XFLPsiam}, we say that 
a generalized inverse $H$ of $A$ is 
\emph{ah-symmetric} (resp. \emph{ha-symmetric})
if $AH$ (resp. $HA$) is symmetric.
It is clear that under the mapping 
$(A,H)\rightarrow(A^\top,H^\top)$, properties \ref{property1} and \ref{property2}  are invariant, while properties \ref{property3} and \ref{property4} are exchanged; so,
although we will concentrate on 
ah-symmetric reflexive generalized inverses,
much of what we work out directly
transfers to the case of 
ha-symmetric reflexive generalized inverses.


We are interested in least-squares problems $\min\{\|A\theta-b\|_2:~\theta\in\mathbb{R}^n\}$,
related to fitting a model of $\beta =\alpha^\top \theta + \epsilon$, where $\epsilon\sim\mathcal{N}(0,\sigma)$. We regard $\alpha$ as a vector of $n$ independent variables, and $\beta$ as the (random) response. We view the rows of $[A,b]$ as $m$ independent realizations of $[\alpha^\top,\beta]$.
We are particularly interested in the situation where  $m>n\gg \rank(A)$.
A very well-known solution of the least-squares problem is $\hat{\theta}:=A^\dag b$, which fits the
model $\hat{\beta}=\alpha^\top \hat{\theta}$.
This solution has a wonderful property --- it quickly gives us the solution for \emph{every}
response vector $b$, and moreover, it is linear in $b$. This all begs the question,
are there other matrices $H$ (besides $A^\dag$) that also have this attractive property? 
It is well known that the answer is yes:  if $H$ is any ah-symmetric generalized inverse, then $\hat{\theta}:=Hb$ solves $\min\{\|A\theta-b\|_2:~\theta\in\mathbb{R}^n\}$ (see \cite[Theorem 6.2.4]{campbell2009generalized} or \cite[Corollary 4.4]{FFL2016}). 
Because a given matrix $A$ does not typically possess a unique ah-symmetric generalized inverse,
there is room to optimize various criteria. 

One important criterion for $H$ is \emph{sparsity}. 
Even if a given matrix is sparse, its M-P pseudoinverse can be completely dense, often leading to a high computational burden in its applications involving many response vectors $b$, especially when we are dealing with a large matrix $A$.  
A sparse $H$ leads to efficiency in calculating $\hat{\theta}:=Hb$.
We can seek such a matrix $H$ with a minimal number of nonzero elements, by applying the standard minimization of its (vector) 1-norm $\|H\|_1:=\sum_{ij} |H_{ij}|$, as a surrogate for the nonconvex ``0-norm''. Besides aiming to induce sparsity, minimization of $\|H\|_1$, or
even a different norm, has the important effect of keeping the elements of $H$ under control, which
reduces numerical errors in the calculation of $\hat{\theta}:=Hb$.

We are further interested in \emph{structured sparsity} for $H$;
specifically, such an $H$ with a \emph{limited number of nonzero rows}.
In evaluating $\hat{\theta}:=Hb$, potentially for many $b$, it is easy to do this very efficiently  
when $H$ has a limited number of nonzero rows,
taking advantage of sparsity without
having to use elaborate sparse data structures. 
Specifically, if $\hat{H}:=H[S,\cdot]$ contains all of the 
nonzeros of $H$, then $\hat{\theta}_{\thickbar{S}}=0$, and $\hat{\theta}_S=\hat{H}b$.
Moreover, for such an $H$, the nonzero rows align with a corresponding 
limited number of columns (i.e., features) of the data matrix $A$, 
and so $\hat{\theta}$ is sparse, and the 
linear model $\hat{\beta}=\alpha^\top \hat{\theta}$ is more \emph{explainable}.

A specific example where we want to evaluate 
$Hb$ for many right-hand sides $b$ is in an ADMM approach to the least-absolute-deviations problem $\min\{\|A\theta-b\|_1:~\theta\in\mathbb{R}^n\}$, 
which aims at a robust fit, 
under the 
assumption that $m>n\gg \rank(A)$;
see \cite[Section 6.1]{boyd2011distributed} (where they assume that $A$ has full column rank; but that assumption is not needed even in their approach, if they simply replace 
$(A^\top A)^{-1}A^\top$ with $A^\dag$). 

There are different ways to work with the criterion of 
wanting $H$ to have a limited number of rows, i.e. ``row-sparsity''.
For example, ($i$) we can try to induce it by minimizing
an appropriate norm, in this case the 2,1-norm  $\|H\|_{2,1}:=\sum_i 
\left\|H_{i\cdot}\right\|_2$\,, or ($ii$) we can impose this structure
and then further seek to (locally) minimize some other criterion, e.g., the (vector) 1-norm   or the 2,1-norm. In both  approaches, we gain an important additional benefit (mentioned above);
by keeping some norm of $H$ under control, we can expect to gain some numerical stability in the calculation of $Hb$.

With the same level of sparsity, structured sparsity should always be preferred to unstructured sparsity. Additionally, algorithms aimed at inducing structured (resp., unstructured) sparsity may or 
may not actually achieve a high level of structured (resp., unstructured) sparsity. Finally, different algorithms with different (but related) goals may
of course have significantly different running times. Therefore,
as we will do in Section \ref{sec:Exper}, we compare the behavior
of different algorithms having different (but related) goals,
on all relevant measures. 

Another important criterion is the \emph{rank} of $H$.
While the sparsity or row-sparsity of $H$ can be viewed as
a kind of simplicity for $H$, such a viewpoint is 
basis dependent. If we make an invertible linear transformation 
of the columns of $A$, that is $A\rightarrow AT$, for 
an invertible and possibly dense $T_{n\times n}$\,, then with an ah-symmetric 
generalized inverse $H$ of $A$, we have $\hat{\theta}=T^{-1}Hb$.
But this matrix $T^{-1}H$ loses the sparsity of $H$. 
So we are also interested in a criterion for simplicity of 
$H$ that is not basis dependent. The natural such criterion is the
\emph{rank} of $H$, which has been employed in many contexts,
in tandem with sparsity criteria
(such as sparse PCA (see \cite{sparcepca}); 
sparse Gaussian mixture models (see \cite{mixtures}, and the references therein); 
low-rank/sparsity matrix decomposition (see \cite{pablo});
low-rank graphical models (see \cite{MohanThesis})), 
to capture a different kind of simplicity/explainability than 
is captured by sparsity. 

In many contexts, rank minimization is naturally induced by minimizing the nuclear norm
(i.e., the 1-norm of the vector of singular values). But in our context,
we have an easier path.
It is well known that: ($i$) if $H$ is a generalized inverse of $A$,
then $\mathrm{rank}(H)\ge\mathrm{rank}(A)$, and ($ii$) a generalized inverse $H$ of $A$ is reflexive if and only if $\mathrm{rank}(H)=\mathrm{rank}(A)$ (see \cite[Theorem 3.14]{RohdeThesis}).
So the minimum possible rank for a generalized inverse of $A$ is $\mathrm{rank}(A)$,
and we can achieve minimization of $\mathrm{rank}(H)$ by simply enforcing \ref{property2},
that is, that $H$ be a reflexive generalized inverse (e.g.,
$A^\dag$ is a reflexive generalized inverse, so it has minimum rank). 
Now, it appears that \ref{property2}
is \emph{nonlinear} in $H$ (i.e., it is a generic system of quadratic 
equations) typically defining a nonconvex region, making enforcing it (in the context of \emph{exact global optimization})
very difficult. But, it is easy to check and well known that 
\ref{property2} becomes linear under properties \ref{property1} and 
\ref{property3}. Specifically,  \ref{property1} and 
\ref{property3} imply that $AH=AA^{\dag}$, and so 
\ref{property2} becomes $HAA^{\dag}=H$.

In summary, our overarching goal is to efficiently calculate 
a sparse or row-sparse reflexive ah-symmetric generalized inverse $H$,
with a controlled value for some norm, in the setting of  $m\! > \! n \! \gg \! \rank(A)$. We develop several methods for 
achieving this goal. 
We note that because we have two objectives (sparsity and low norm), there is naturally
a trade-off to look at.

\smallskip 
\noindent {\bf Notation.} 
$\|\cdot\|_p$ denotes vector $p$-norm ($p\!\geq\! 1$), $\|\cdot\|_0$ denotes vector ``0-norm''. For matrices,
$\tr(\cdot)$ denotes trace, $\det(\cdot)$ denotes determinant, $\rank(\cdot)$ denotes rank,  $\|\cdot\|_F$ denotes the Frobenius norm, $\mathcal{R}(\cdot)$ denotes the column space.
$A[S,T]$ is the submatrix of $A$ having row indices $S$ and column indices $T$, 
with $S=\cdot$ (resp. $T=\cdot$) indicating 
all rows (resp. columns); further,
$A_{i\cdot}:=A[\{i\},\cdot]$,
$A_{\cdot j}:=A[\cdot,\{j\}]$, and $A_{ij}:=A[\{i\},\{j\}]$. For $A\in\mathbb{R}^{m\times n}$, $\|A\|_{p,q}:=\|(\|A_{1\cdot}\|_p,\ldots,\|A_{m\cdot}\|_p)\|_q$ ($p,q\geq 0$).
$\mathbf{e}_i$ denotes the $i$-th standard unit vector.
$I$ denotes an identity matrix (sometimes  with a subscript indicating its order). For $A\in\mathbb{R}^{m\times n}$, 
for the \emph{compact singular-value decomposition}  (compact SVD), we  write $A=U\Sigma V^\top$,
where $r:=\rank(A)$,  $U_{m\times r}$ and $V_{n\times r}$ are real orthogonal matrices, 
 $\Sigma_{r\times r}$ is a diagonal matrix, with positive diagonal entries
 $\sigma_1\geq\sigma_2\geq \cdots \geq \sigma_r$\,, known as the 
 singular values of $A$. For the  \emph{(full) singular-value decomposition}  (SVD),
we instead have $U_{m\times m}$\,, $V_{n\times n}$\,, 
 with nonnegative diagonal entries
 $\sigma_1\geq\sigma_2\geq \cdots \geq \sigma_{\min\{m,n\}}$\,. We denote by $\argmin\{\cdot\}$
\emph{any} solution of the associated minimization problem.

\smallskip 
\noindent {\bf Literature.}
For basic information on generalized inverses, we refer to  \cite{RohdeThesis}, \cite{BenIsrael1974}, and \cite{campbell2009generalized}.
%
%
\cite{dokmanic} (and then \cite{dokmanic1,dokmanic2})
introduced the idea of seeking to induce sparsity in $H$ by linear optimization; specifically, they minimize the (vector) 1-norm  over the set of left inverses (i.e., $HA=I_n$) or right inverses
(i.e., $AH=I_m$).  \cite{dokmanic1} also considered the idea of inducing row-sparsity by minimizing the 2,1-norm over the same sets of matrices. Additionally, \cite{dokmanic1} presents ADMM algorithms for these minimization  problems, using specific projections for the cases where the matrices are full row- or column-rank. 
It is important to realize that the left (resp., right) inverse case applies
only to full column (resp., row) rank matrices $A$. 
In particular, every left inverse $H$ is a reflexive ha-symmetric generalized inverse, 
and so for such an $H$,
$\hat{\theta}:=Hb$ solves the
least-squares problem (for all $b$) \emph{only} when $H=A^\dag$.
In fact, we are motivated by the case in which $A$ is neither full row rank nor full column rank, and so
the left/right inverse approach does not apply at all in our setting.

In \cite{FFL2016}, we 
proposed the idea of
finding various types of sparse generalized inverses
by imposing combinations of the M-P properties in the context of straight-forward linear-optimization formulations,
minimizing the (vector) 1-norm of $H$.
But we did not propose any practical idea for
controlling the rank of solutions.
Moreover, the approaches proposed did not scale well. 

In \cite{FampaLee2018ORL}, we
introduced the idea of 
imposing a block structure on $H$, 
which implicitly enforces \ref{property1} and \ref{property2},
and then carrying out a combinatorial local search, giving an approximation guarantee on (vector) 1-norm minimization for reflexive generalized inverses. 
That work does not directly related to 
our present setting because it ignores 
property \ref{property3}.
Additionally, at the time, the empirical quality (1-norm) of the solutions
for that method on large instances was unknown, because we could not solve 
the large linear-optimization formulations
that would give us lower bounds.

In \cite{ojmoFLP21},
we  proposed
different algorithmic approaches  that start from a (vector) 1-norm minimizing ah-symmetric generalized inverse, and gradually decrease its rank, by iteratively imposing the reflexive property. The algorithms iterate
 until the ah-symmetric generalized inverse has the least possible rank producing intermediate solutions during the iterations, trading
 off low 1-norm against low rank. The best approach investigated was a cutting-plane method that solves linear-optimization problems at each iteration, although it is capable of  constructing interesting intermediate solutions (trading off rank and 1-norm), it does  not  scale well, as the solution of many dense linear-optimization problems  is required. Only square matrices with size $50$ and  rank $25$  were considered in the numerical experiments.

In \cite{XFLPsiam}, we significantly extended the work in \cite{FampaLee2018ORL}, 
in particular to approximate 
(vector) 1-norm minimizing ah-symmetric reflexive
generalized inverses, using a 
column-block construction, which implicitly enforces \ref{property1},
\ref{property2},
and \ref{property3}, and
maintains structured row-sparsity during a local search. 
Comparisons between solutions of the local-search procedure and the minimal 1-norm of ah-symmetric generalized inverses were presented  for matrices of size up to $120\times 60$ and rank $30$. It was also mentioned that, within a time limit of two hours, the solutions of the linear-optimization problems to minimize the 1-norm of generalized  and ah-symmetric generalized inverses were obtained for only one out of five given  matrices of size $200\times 100$ and rank $50$. The results indicated that, in practice, the solutions of the linear-optimization problems were not easy for large matrices. \cite{XFLrank12} is a companion 
work that in particular analyzes the motivating cases for which $\rank(A)\in\{1,2\}$; while the case of rank one is trivial, the complexity of the case of rank two reveals that there is no simple solution in general.   

In \cite{FLPXjogo}, we carried out a detailed computational study of local-search procedures based on the results of \cite{FampaLee2018ORL} and  \cite{XFLPsiam}. An experimental analysis of the procedures was performed, but comparisons between the local search solutions and   the minimum 1-norm of generalized  and ah-symmetric 
generalized inverses were only presented for square matrices of size up to $100$ and rank $50$, again because of the intractability of the large dense linear-optimization problems. 

In \cite{PFLX_ORL}, in addition to considering 1-norm minimization with linear optimization to induce (unstructured) sparsity of generalized inverses, which was first advanced by \cite{FFL2016}, we also considered 2,1-norm minimization with second-order-cone optimization to induce row-sparsity.  Furthermore, \cite{PFLX_ORL} showed how to make the ideas of linear and second-order-cone optimization much more efficient/scalable. A key  idea exploited  in  \cite{PFLX_ORL} 
is as follows.
Let $A \in \mathbb{R}^{m \times n}$ with rank $r$ and $A =: U \Sigma V^\top$ be the (full) SVD of $A$, where $U \in \mathbb{R}^{m \times m}$, $V \in \mathbb{R}^{n \times n}$ are orthogonal matrices  ($U^\top U = I_m, V^\top V= I_n$), and $\Sigma \in \mathbb{R}^{m \times n}$ with 
\[
\Sigma =: \begin{bmatrix}\underset{\scriptscriptstyle r\times r}{D} & \underset{\scriptscriptstyle r\times (n-r)}{0}\\
\underset{\scriptscriptstyle (m-r)\times r}{0} & \underset{\scriptscriptstyle (m-r)\times (n-r)}{0}\end{bmatrix},
\]
\noindent with $D$ being a diagonal matrix with rank $r$. 
Let $H \in \mathbb{R}^{n \times m}$ and $\Gamma:=V^\top H U$, where
\[
\Gamma=: \begin{bmatrix}\underset{\scriptscriptstyle r\times r}{X} & \underset{\scriptscriptstyle r\times (m-r)}{Y}\\
\underset{\scriptscriptstyle (n-r)\times r}{Z} & \underset{\scriptscriptstyle (n-r)\times (m-r)}{W}\end{bmatrix},
\]
then $H=I_n H I_m = (VV^\top) H (U U^\top) = V\Gamma U^\top$.

\vbox{
\begin{lemm}[\protect{see \cite{BenIsrael1974}, p.\! 208, \!Ex. \!\!14 (no~proof); see \cite[Lemmas 2-5]{PFLX_ORL} (with proof)}]\label{lem:structure}\phantom{a}
\begin{itemize}
\item 
 \ref{property1} is equivalent to  
 $X = D^{-1}$.
\item
  If \ref{property1} is satisfied, then \ref{property2} is equivalent to 
  $ZDY = W$.
\item 
If \ref{property1} is satisfied, then \ref{property3} is equivalent to 
    $Y = 0$.
\item
    If \ref{property1} is satisfied, then \ref{property4} is equivalent to  
    $Z = 0$.
\end{itemize}
\end{lemm}
}

So, for the purpose of enforcing any subset of the M-P properties on a generalized inverse $H$, except the case of \ref{property2} alone, we simply
set appropriate blocks of $\Gamma$ to zeros. 

Consider the following natural convex-optimization problems, aimed at inducing sparsity and row-sparsity for ah-symmetric reflexive generalized inverses, respectively.
\begin{align*}
\tag{\mbox{$P^1_{123}$}}
\label{prob:min1norm}
\min_{H\in\mathbb{R}^{n \times m}} \left\{\|H\|_{1}
~:~ \ref{property1},\ref{property2},\ref{property3}
\right\}
\end{align*}
\noindent and 
\begin{align*}
\tag{\mbox{$P_1^{2,1}$}}
\label{prob:min21norm}
\min_{H\in\mathbb{R}^{n \times m}}\left\{\|H\|_{2,1} ~:~
\ref{property1}\right\}.
\end{align*}
It is important to note that due to \cite[Corollary 9]{PFLX_ORL}, optimal solutions of \ref{prob:min21norm}
automatically satisfy 
\ref{property2} and 
\ref{property3} (i.e., they are reflexive and ah-symmetric),
which we exploit in establishing Theorem \ref{thm:reduced} below.

 Let 
\[
V:= \begin{bmatrix}\underset{\scriptscriptstyle n\times r}{V_1} & \underset{\scriptscriptstyle n\times (n-r)}{V_2}\end{bmatrix}, ~ U:= \begin{bmatrix}\underset{\scriptscriptstyle m\times r}{U_1} & \underset{\scriptscriptstyle m\times (m-r)}{U_2}\end{bmatrix}. 
\]
We remark that 
    $U_1$ and $V_2$ are full column rank matrices and $U_1^\top U_1 = I$ and $V_2^\top V_2 = I$, then $U_1^\dagger = (U_1^\top U_1)^{-1}U_1^\top = U_1^\top$ and $V_2^\dagger = (V_2^\top V_2)^{-1} V_2^\top = V_2^\top$\,.

From Lemma \ref{lem:structure}, we can deduce the following result. 
\begin{theo}[\protect{\cite[Section 5]{PFLX_ORL}}]\label{thm:reduced}
    \ref{prob:min1norm} and \ref{prob:min21norm}, respectively, 
can be efficiently reformulated  as 
\begin{equation}\tag{\mbox{$\mathcal{P}^1_{123}$}}\label{prob:barmin1norm}
z(\mathcal{P}^1_{123}):=\displaystyle\min_{Z\in\mathbb{R}^{(n-r) \times r}}\left\| V_1D^{-1}U_1^\top + V_2ZU_1^\top\right\|_1\, ,
\end{equation}
\begin{equation}\tag{\mbox{$\mathcal{P}_1^{2,1}$}}\label{prob:barmin21norm}
z(\mathcal{P}_1^{2,1}):=\displaystyle\min_{Z\in\mathbb{R}^{(n-r) \times r}} 
\left\| V_1 D^{-1} + V_2 Z \right\|_{2,1}\, .
\end{equation}
\end{theo}
\smallskip

The more compact formulations presented in Theorem \ref{thm:reduced} were exploited in the numerical experiments performed in \cite{PFLX_ORL} to increase the size of matrices for which the solutions of \ref{prob:min1norm} and \ref{prob:min21norm} were known. Within a time limit of 5 hours, the solutions of \ref{prob:barmin1norm} and \ref{prob:barmin21norm} were obtained for instances with $(m,n,r)$ equal to $(280,140,60)$ and $(2000,1000,500)$, respectively. These results allowed us to finally evaluate the empirical quality of the proposed local-search solutions on larger instances.

Concerning the contribution of our present work, we  take a further step in  efficiently obtaining solutions to \ref{prob:min1norm} and \ref{prob:min21norm}. We leverage  Theorem \ref{thm:reduced} once again  and  develop very efficient/scalable ADMM (Alternating Direction Method of Multipliers) algorithms to solve them.  Furthermore, we  develop two ADMM algorithms aiming at imposing structured sparsity by limiting the number of nonzero rows for an ah-symmetric 
reflexive generalized inverse. For the first,  we do \emph{not} consider the minimization of a norm of the nonzero submatrix of $H$, leading to a fast ADMM algorithm that can take advantage of the efficient solutions presented in  the literature for the subproblems solved. For the second, we minimize the 2,1-norm of $H$, in addition to limiting its number of nonzero rows. We provide a closed-form solution for the nonconvex subproblem solved, deriving an efficient ADMM algorithm for this problem as well. Although there is no guarantee of the convergence of these two ADMMs to solutions of the nonconvex problems addressed, we demonstrate their efficacy through our numerical results.    
We show that: (i) our ADMM approaches are much more scalable than direct methods
aimed at 
\ref{prob:barmin1norm}
and \ref{prob:barmin21norm}
(using general-purpose software like \texttt{Gurobi} and \texttt{MOSEK}), and (ii) our ADMM methods provide a good complement to our local-search methods. 
They not only provide a means to better evaluate the solutions computed with the local-search methods, but can in fact construct solutions of comparable quality with respect to sparsity and row-sparsity in less time.
Finally,  we also demonstrate in this work that the
efficient minimum 
(vector) 1-norm $r$-approximating  local search
of \cite{XFLPsiam},
for an ah-symmetric 
generalized inverse,
in fact  gives a factor-$r$
approximation for  
2,1-norm minimization as well, which aims directly at inducing row sparsity.

\bigskip 
\noindent {\bf Organization.} 
In \S\ref{sec:ADMM}, we  develop the 
 ADMM 
 algorithms aimed at inducing  unstructured and structured sparsity in our setting. 
In \S\ref{sec:LS}, we give an approximation guarantee on 2,1-norm 
minimization for ah-symmetric 
generalized inverses, for  the
 minimum 
(vector) 1-norm $r$-approximating  local search
of \cite{XFLPsiam}. For very low $r$, 
we do not need to settle for approximation, and 
so we make a detailed analysis of the cases of rank one and two in 
Appendix \ref{app:oneandtwo}. 
In \S\ref{sec:admm20}, we develop the two ADMM algorithms aiming at imposing structured sparsity by limiting the number of nonzero rows for an ah-symmetric 
reflexive generalized inverse.     
In \S\ref{sec:Exper}, we present the results of computational experiments, demonstrating the favorable performance of our new methods. 
 In \S\ref{sec:Next}, we indicate some next steps for this line of research.


\section{Inducing sparsity and row sparsity with ADMM}\label{sec:ADMM}

It is natural to consider 
specialized algorithms for attacking problems like 
\ref{prob:min1norm} and
\ref{prob:min21norm}, seeking fast convergence to near optima. It is even more enticing to 
seek such methods for their compact forms, 
\ref{prob:barmin1norm}
and \ref{prob:barmin21norm}. 
In these compact forms, we have unconstrained minimization problems in
a single variable.
Motivated by for example 
\cite[Section 6.1]{boyd2011distributed},
we can introduce a second variable and linear linking constraints
to  seek to develop efficient ADMM algorithms. 

\subsection{ADMM for 1-norm minimization}\label{subsec:admm1}
We seek to develop an ADMM algorithm for \ref{prob:barmin1norm}\,. Initially, by introducing a variable $E \in \mathbb{R}^{n \times m}$, we rewrite  \ref{prob:barmin1norm}  as
\begin{equation}\label{prob:admm_mat_1norm}
\min \left\{
 \|E\|_{1} ~:~  E =  V_1D^{-1}U_1^\top + V_2ZU_1^\top
 \right\}.
\end{equation}
The augmented Lagrangian function associated to \eqref{prob:admm_mat_1norm} is
\begin{align*}
\mathcal{L}_\rho(Z,E,\Lambda)\!
 &:=\! 
 \|E\|_{1} \!+  \!\langle \Theta, V_1D^{-1}U_1^\top \!+\! V_2ZU_1^\top \!-\!E \rangle+\! \textstyle \frac{\rho}{2}\!\left\|V_1D^{-1}U_1^\top \!+\! V_2ZU_1^\top \!-\!E \right\|^2_F\\
&\phantom{:}=\! 
 \|E\|_{1} \!+\! \textstyle \frac{\rho}{2}\!\left\|V_1D^{-1}U_1^\top \!+\! V_2ZU_1^\top \!-\!E \!+\!\Lambda \right\|^2_F - 
  \frac{\rho}{2}\left\|\Lambda\right\|_F^2\,,
\end{align*}
where $\rho >0$ is the penalty parameter,
$\Theta\in \mathbb{R}^{n \times m}$ is the Lagrangian multiplier,
and $\Lambda $  is the scaled  Lagrangian multiplier, that is  $\Lambda := \Theta/\rho$. We will apply the ADMM method to \ref{prob:barmin1norm}\,, by iteratively solving, for $k=0,1,\ldots,$ 
\begin{align}
    &Z^{k+1}:=\textstyle\argmin_Z ~ \mathcal{L}_\rho(Z,E^{k},\Lambda^k),\label{eq:Zmin1normsubpa}\\
    &E^{k+1}:=\textstyle\argmin_E ~ \mathcal{L}_\rho(Z^{k+1},E,\Lambda^k),\label{eq:Emin1normsubpa}\\
    &\textstyle\Lambda^{k+1}:=\Lambda^{k} + V_1D^{-1}U_1^\top + V_2 Z^{k+1}U_1^\top - E^{k+1}.\nonumber
\end{align}
Next, we detail how to solve the subproblems above.

\medskip

\noindent {\bf{Update \texorpdfstring{$Z$}{Z}:}} 
To update $Z$,  we consider  subproblem \eqref{eq:Zmin1normsubpa}, more specifically,
\begin{align}\label{eq:Zmin1normsubprob}
Z^{k+1}:=\textstyle \argmin_{Z}\{  \left\| J - V_2 Z U_1^\top\right\|^2_F\},
\end{align}
where $J:=  E^k -V_1D^{-1}U_1^\top -\Lambda^k$. We can easily  verify that the solution of \eqref{eq:Zmin1normsubprob} is given by $Z^{k+1} =  V_2^\top J U_1$\,.


\medskip

\noindent {\bf{Update \texorpdfstring{$E$}{E}:}} 
To update $E$,   we consider  subproblem \eqref{eq:Emin1normsubpa}, more specifically, 
\begin{align}
E^{k+1}:=\textstyle\argmin_E\{ \|E\|_{1} + \frac{\rho}{2}\left\|E -Y  \right\|^2_F\},\label{eq:Emin1normsubprob}
\end{align}
where $Y := V_1D^{-1}U_1^\top + V_2Z^{k+1}U_1^\top  +\Lambda^k$.

\begin{prop}
[\protect{see \cite[Section 4.4.3]
{boyd2011distributed}}]\label{lem:closedformula1norm}
The solution of \eqref{eq:Emin1normsubprob} is given by 
$
    E^{k+1}_{ij} = S_{1\!/\!\rho}(Y_{ij}),
    $
for $i = 1,\dots,n$ and $j = 1,\dots,m$, where  the soft thresholding operator $S$ is defined as
\begin{align*}
    S_\kappa(a) := \begin{cases}
    a - \kappa, \quad & a > \kappa;\\ 
    0, & |a| \leq \kappa;\\ 
    a + \kappa, &a < -\kappa.
\end{cases}
\end{align*}
\end{prop}

\medskip

\noindent
{\bf{Initialization of the variables:}}   
We need to initialize $\Lambda$ and $E$. First, we set $\Lambda^0:=\hat\Theta/\rho$, where $\textstyle \hat \Theta :=\frac{1}{\|V_1 U_1^\top\|_{\infty}} V_1 U_1^\top$. Our goal was to select the Lagrangian multiplier   $\hat\Theta$ as an easily computable feasible solution to  the dual problem of  \eqref{prob:admm_mat_1norm}, given by
\[
\max_{ \Theta}\{\tr(D^{-1}V_1^\top \Theta U_1):V_2^\top \Theta U_1 = 0,~\|\Theta\|_{\infty} \leq 1\}.
\]
Due to the orthogonality of $V$, we can easily verify the feasibility of $\hat\Theta$.
Furthermore, it gives a positive objective value and has maximum infinity norm;
in contrast to the zero matrix which is also feasible.

Then, from Lemma \ref{lem:structure}, we note that the ah-symmetric reflexive generalized inverses of $A$ can be written as $V_1D^{-1}U_1^\top + V_2ZU_1^\top$, and  the M-P pseudoinverse $A^\dagger$ can be written as $V_1D^{-1}U_1^\top$. We recall that $A^\dagger$ is the generalized inverse of $A$ with minimum Frobenius norm. Then, aiming to obtain $Z^1=0$ when solving \eqref{eq:Zmin1normsubprob} in the first iteration of the algorithm, and consequently starting the algorithm with a Frobenius norm minimizing ah-symmetric reflexive generalized inverse, we set $E^0 := V_1D^{-1}U_1^\top + \Lambda^0$.

\medskip

\noindent{\bf{Stopping criterion:}} 
We consider a stopping criterion  from \cite[Section 3.3.1]{boyd2011distributed}, and select an absolute tolerance $\epsilon^{\mbox{\scriptsize abs}}$ and a
relative tolerance $\epsilon^{\mbox{\scriptsize rel}}$. The algorithm stops at iteration $k$ if  
\begin{align}
    &\|r^k\|_F \leq \epsilon^{\mbox{\scriptsize abs}}\sqrt{nm} 
    + \epsilon^{\mbox{\scriptsize rel}}\max\left\{\|E^k\|_F\,, \|V_2Z^kU_1^\top\|_F\,,\|V_1D^{-1}U_1^\top\|_F\right\}, \label{sc1r} \\
    &\|s^k\|_F \leq \epsilon^{\mbox{\scriptsize abs}}
    \sqrt{(n-r)r} +    
\epsilon^{\mbox{\scriptsize rel}}\rho \|V_2^\top\Lambda^kU_1\|_F\,,\label{sc1s}
\end{align}
where $r^k := V_1D^{-1}U_1^\top + V_2Z^kU_1^\top - E^k$ is the primal residual,  and $s^k:= \rho V_2^\top(E^{k} - E^{k-1})U_1$\, is the dual residual.

\bigskip

In Algorithm \ref{alg:admm1norm}, we present the ADMM algorithm for \ref{prob:barmin1norm}\,. We observe that in steps 7 to 9, the elements of $E^{k+1}$ can be computed in parallel.

\begin{algorithm2e}[!ht]
\footnotesize{
\KwIn{ $A\in \mathbb{R}^{m\times n}$,
$\Lambda^0 \in \mathbb{R}^{n\times m}$, $E^0 \in \mathbb{R}^{n\times m}$, $\rho > 0$.}
\KwOut{$H\in \mathbb{R}^{n\times m}$.}
$U,\Sigma,V := \texttt{svd}(A)$, $k:=0$\;
Get $U_1,V_1,V_2,D^{-1}$ from $U,\Sigma,V$\;
\While {not converged} 
{
$J := E^k -V_1D^{-1}U_1^\top - \Lambda^k$\;
$Z^{k+1} := V_2^\top JU_1$\;
$Y := V_1D^{-1}U_1^\top + V_2 Z^{k+1}U_1^\top + \Lambda^k$\;
\For{$i = 1,\dots,n$}{
\For{$j = 1,\dots,m$}{
$E^{k+1}_{ij} :=  S_{1\!/\!\rho}(Y_{ij})$; \qquad (see Proposition \ref{lem:closedformula1norm}){\color{white}\;}
}
}
$\Lambda^{k+1}:=\Lambda^{k} + V_1D^{-1}U_1^\top + V_2 Z^{k+1}U_1^\top - E^{k+1}$\;
$k:=k+1$\;
}
$H:=V_1D^{-1}U_1^\top + V_2Z^kU_1^\top$\;
\caption{ADMM for \ref{prob:barmin1norm} (ADMM$_{1}$)}\label{alg:admm1norm}
\hypertarget{algadmm1norm}{}
}
\end{algorithm2e}

\subsection{ADMM for 2,1-norm minimization}\label{subsec:ADMM21}

Next, we apply ADMM for solving \ref{prob:barmin21norm}. Initially, by introducing a variable
$E \in \mathbb{R}^{n \times r}$, we rewrite  \ref{prob:barmin21norm}  as
\begin{equation}\label{prob:admm_mat_21norm}
\min
\left\{\|E\|_{2,1} ~:~  E =  V_1D^{-1} + V_2 Z\right\}.
\end{equation}

\begin{rema}
    In Theorem \ref{thm:reduced}, we use the orthogonality of $U_1$ to replace $\|V_1D^{-1}U_1^\top + V_2 ZU_1^\top\|_{2,1}$ with $\|V_1D^{-1} + V_2 Z\|_{2,1}$ in the objective function of \ref{prob:barmin21norm}. This leads to a very helpful reduction in the dimension of the variable $E$ in the ADMM algorithm presented in this section, compared to what we presented in the previous section.
\end{rema}

\medskip

The augmented Lagrangian function associated to \eqref{prob:admm_mat_21norm} is
\begin{align*}
&\mathcal{L}_\rho(Z,E,\Lambda)\!:=\! \textstyle
\|E\|_{2,1} \!+ \!\frac{\rho}{2}\left\|V_1D^{-1} \!+\! V_2Z \!-\!E \!+\! \Lambda \right\|^2_F  - 
  \frac{\rho}{2}\left\|\Lambda\right\|_F^2\,,
\end{align*}
where $\rho >0$ is the penalty parameter and $\Lambda \in \mathbb{R}^{n \times r}$ is the scaled Lagrangian multiplier.  We will apply the ADMM method to \ref{prob:barmin21norm}, by iteratively solving, for $k=0,1,\ldots,$  
\begin{align}
    &Z^{k+1}:=\textstyle\argmin_Z ~ \mathcal{L}_\rho(Z,E^{k},\Lambda^k),\label{eq:Zmin21normsubpa}\\
    &E^{k+1}:=\textstyle\argmin_E ~ \mathcal{L}_\rho(Z^{k+1},E,\Lambda^k),\label{eq:Emin21normsubpa}\\
    &\textstyle\Lambda^{k+1}:=\Lambda^{k} + V_1D^{-1} + V_2 Z^{k+1} - E^{k+1}\,.\nonumber
\end{align}
Next, we detail how to solve the subproblems above.

\medskip

\noindent{\bf{Update of \texorpdfstring{$Z$}{Z}:}} 
To update $Z$,  we consider  subproblem \eqref{eq:Zmin21normsubpa}, that is
\begin{equation}\label{eq:Zmin21normsubprob}
Z^{k+1}:=\textstyle\argmin_Z  \left\| J - V_2 Z\right\|^2_F\,,
\end{equation}
where $J:=  E^k -V_1D^{-1} -\Lambda^k$. Similarly to \eqref{eq:Zmin1normsubprob}, we can verify that $Z^{k+1} = V_2^\top J$ is an optimal  solution to  \eqref{eq:Zmin21normsubprob}.

\medskip

\noindent{\bf{Update of \texorpdfstring{$E$}{E}:}} 
To update $E$,   we consider  subproblem \eqref{eq:Emin21normsubpa}, that is 
\begin{align}
E^{k+1}:=\textstyle\argmin_E\{ \|E\|_{2,1} + \frac{\rho}{2}\left\|E -Y  \right\|^2_F\},\label{eq:Emin21normsubprob}
\end{align}
where $Y := V_1D^{-1} + V_2Z^{k+1}  +\Lambda^k$\,.

\begin{prop}[\protect{see \cite[Proposition 1]{yuan2006model}}]\label{lem:closedformula21norm} The solution of \eqref{eq:Emin21normsubprob} is given by 
\begin{equation}\label{updateE}
    E^{k+1}_{i\cdot} := \begin{cases}
    \frac{ \|Y_{i\cdot}\|_2 - 1\!/\!\rho}{\|Y_{i\cdot}\|_2} Y_{i\cdot}\,,\quad &\text{if } 1\!/\!\rho < \|Y_{i\cdot}\|_2 \,;
    \\
    0,&\text{otherwise.}
\end{cases}
\end{equation}
\end{prop}

\medskip

\noindent{\bf{Initialization of the variables:}} 
To initialize the variables, we apply the same ideas discussed in \S\ref{subsec:admm1}. Considering now the dual problem of \eqref{prob:admm_mat_21norm},
\begin{equation*}
\max\limits_{ \Theta}\{ \tr(D^{-1}V_1^\top \Theta):V_2^\top \Theta = 0,\|\Theta_{i\cdot}\|_{2} \!\leq\! 1,i\! =\! 1,\dots,n\},
\end{equation*}
we set  $\hat \Theta
:= V_1/ \kappa$, where $\kappa := \displaystyle \max_{i = 1,\dots,n}\|{V_1[i,\cdot]}\|_{2}$\,,   $\Lambda^0 :=\hat\Theta/ \rho$\,,  and $E^0 := V_1D^{-1} + \Lambda^{0}$. 

\medskip

\noindent{\bf{Stopping criterion:}} 
We adopt the same stopping criterion described in \S\ref{subsec:admm1}. With the primal and dual residuals at iteration $k$ given respectively by 
$r^k := V_1D^{-1} + V_2Z^k - E^k$ and $s^k:= \rho V_2^\top(E^{k} - E^{k-1})\,$, the criterion becomes
\begin{align}
    &\|r^k\|_F \leq \epsilon^{\mbox{\scriptsize abs}}\sqrt{nr} 
    + \epsilon^{\mbox{\scriptsize rel}}\max\left\{\|E^k\|_F\,, \|V_2Z^k\|_F\,,\|V_1D^{-1}\|_F\right\},\label{sc2r}\\
    &\|s^k\|_F \leq \epsilon^{\mbox{\scriptsize abs}}
    \sqrt{(n-r)r}  
    + \epsilon^{\mbox{\scriptsize rel}}\rho \|V_2^\top\Lambda^k\|_F\,.\label{sc2s}
\end{align}

\bigskip

In Algorithm \ref{alg:admm21norm}, we present the ADMM algorithm for \ref{prob:barmin21norm}\,.

\begin{algorithm2e}[!ht]
\footnotesize{
\KwIn{$A\in \mathbb{R}^{m\times n}$,
$\Lambda^0 \in \mathbb{R}^{n\times r}$, $E^0 \in \mathbb{R}^{n\times r}$, 
  $\rho > 0$.}
\KwOut{$H\in \mathbb{R}^{n\times m}$.}
$U,\Sigma,V := \texttt{svd}(A)$, $k:=0$\;
Get $U_1,V_1,V_2,D^{-1}$ from $U,\Sigma,V$\;
\While {not converged} 
{
$J := E^k -V_1D^{-1} - \Lambda^k$\;
$Z^{k+1} := V_2^\top J$\;
$Y := V_1D^{-1} + V_2 Z^{k+1} + \Lambda^k$\;
\For{$i = 1,\dots,n$}{
    \bf{if} $\|Y_{i\cdot}\|_2> 1\!/\!\rho$ 
   \bf{then} $E^{k+1}_{i\cdot} := \frac{\|Y_{i\cdot}\|_2 - 1\!/\!\rho}{\|Y_{i\cdot}\|_2} Y_{i\cdot}$ \!; ~ (see Proposition \ref{lem:closedformula21norm}){\color{white}\;}
    \bf{else}
    $E^{k+1}_{i\cdot} := 0$\;
}
 $\Lambda^{k+1} := \Lambda^k + V_1D^{-1} + V_2 Z^{k+1} - E^{k+1}$\;
 $k:=k+1$\;
}
$H:=V_1D^{-1}U_1^\top + V_2Z^kU_1^\top$\;
\caption{ADMM for \ref{prob:barmin21norm} (ADMM$_{2,1}$)\label{alg:admm21norm}}
\hypertarget{algadmm21norm}{}
}
\end{algorithm2e}


\section{Imposing row sparsity with column-block solutions and local search}
\label{sec:LS}

\cite[Section 3]{XFLPsiam} developed a local-search based approximation algorithm to efficiently calculate an approximate (vector) $1$-norm minimizing ah-symmetric generalized inverse with $r:=\rank(A)$ nonzero rows. 
It is very natural to investigate how well this approximation algorithm does for other norms, in particular the 2,1-norm, which we emphasize in our present work. 
Here we show that the same matrix calculated by this approximation algorithm is also a good approximate 2,1-norm minimizing ah-symmetric generalized inverse. We note that our approximation factor for the 2,1-norm is much better than what would be obtained by simply applying norm inequalities to the result for the 1-norm (see Remark \ref{rem:approx_from1_to_21}).
Our local search works with  ``column-block solutions'', as presented in the following result. 
\begin{theo}[\protect{\cite[Theorem 3.3]{XFLPsiam}}]\label{thm:ahconstruction}
For $A\!\in\!\mathbb{R}^{m\times n}$, let $r \!:= \!\rank(A)$. For any $T$, an ordered subset of $r$ elements from $\{1,\dots,n\}$, let $\hat{A}\!:=\!A[\cdot,T]$ be the $m \times r$ submatrix of $A$ formed by columns $T$. If $\rank(\hat{A})\!=\!r$, let
$
\hat{H}\! := \!\hat{A}^{\dagger}\! \!=\!(\hat{A}^\top\hat{A})^{-1}\hat{A}^\top.
$
Then the $n \times m$ matrix $H$ with all rows equal to zero, except rows $T$, which are given by $\hat{H}$, is an ah-symmetric reflexive generalized inverse of $A$.
\end{theo}

In Appendix \ref{app:oneandtwo}, we  consider the cases in which $\rank(A)\in\{1,2\}$. As was done for the 1-norm in \cite{XFLrank12}, we  use the result in Theorem \ref{thm:ahconstruction} to address the exact 2,1-norm minimization of ah-symmetric reflexive generalized inverses for these simpler cases. After observing the difficulty of the exact  minimization already for rank $2$, 
we consider an approximation algorithm for general matrices
of constant rank $r$, as discussed next.

\begin{defi}[\protect{\cite[Definition 3.5]{XFLPsiam}}]\label{def:localsearch_ahsym}
\!\!\! For $A\!\in\!\mathbb{R}^{m\times n}$, let $r \!:= \!\rank(A)$, and~$S$ be an ordered subset of $r$ elements from $\{1,\dots,m\}$ such that $A[S,\cdot]$ has linearly independent rows. For $T$ an ordered subset of $r$ elements from $\{1,\dots,n\}$, and fixed $\epsilon\!\ge \!0$, if $|\det(A[S,T])|$ cannot be increased by a factor of more than $1\!+\!\epsilon$ by swapping an element of $T$ with one from its complement, we say that $A[S,T]$ is a $(1\!+\!\epsilon)$-local maximizer for the absolute determinant on the set of $r\!\times\! r$ nonsingular submatrices of $A[S,\cdot]$.
\end{defi}

\begin{rema}\label{rem:ptime}
 \cite[Theorem 3.9]{XFLPsiam} showed that a local optimum  satisfying Definition \ref{def:localsearch_ahsym}
can be calculated in polynomial time (when $A$ is rational), for any fixed $\epsilon\!>\!0$.
\end{rema}

\medskip

\begin{theo}[\protect{\cite[Theorem 3.7]{XFLPsiam}}]\label{thmrwithP3}
For $A\!\in\!\mathbb{R}^{m\times n}$, let $r \!:= \!\rank(A)$, and let $S$ be an ordered subset of $r$ elements from $\{1,\dots,m\}$ such that $A[S,\cdot]$ has linearly independent rows. Choose $\epsilon\ge0$, and let $\tilde{A}:=A[S,T]$ be a $(1+\epsilon)$-local maximizer for the absolute determinant on the set of $r\times r$ nonsingular submatrices of $A[S,\cdot]$. Then the $n \times m$ matrix $H$ constructed by Theorem \ref{thm:ahconstruction} over $\hat{A}:=A[\cdot,T]$, is an ah-symmetric reflexive generalized inverse of $A$ satisfying $\|H\|_1\le r(1+\epsilon)
 \|H_{opt}^{1}\|_1$\,, where 
$H_{opt}^1$ is a (vector) $1$-norm minimizing ah-symmetric
reflexive generalized inverse of $A$.
\end{theo}

Next, we will show that the constructed $H$ of 
Theorem \ref{thmrwithP3} has 2,1-norm within a factor of $r$ of the 2,1-norm minimizing solution as well. 

\begin{theo} \label{thm:21rwithP3}
For $A\!\in\!\mathbb{R}^{m\times n}$, let $r \!:= \!\rank(A)$, and let $S$ be an ordered subset of $r$ elements from $\{1,\dots,m\}$ such that $A[S,\cdot]$ has linearly independent rows. Choose $\epsilon\ge0$, and let $\tilde{A}:=A[S,T]$ be a $(1+\epsilon)$-local maximizer for the absolute determinant on the set of $r\times r$ nonsingular submatrices of $A[S,\cdot]$. Then the $n \times m$ matrix $H$ constructed by Theorem \ref{thm:ahconstruction}
over $\hat{A}:=A[\cdot,T]$, is an ah-symmetric reflexive generalized inverse of $A$ satisfying $\|H\|_{2,1}\le r(1+\epsilon)\|H_{opt}^{2,1}\|_{2,1}$\,, where $H_{opt}^{2,1}$ is a 2,1-norm minimizing ah-symmetric reflexive generalized inverse of $A$.
\end{theo}

We note that $H_{opt}^{2,1}$ defined in Theorem \ref{thm:21rwithP3} is also a 2,1-norm minimizing generalized inverse of $A$ because as stated in \S\ref{sec:introd}, optimal solutions of \ref{prob:min21norm}
automatically satisfy 
\ref{property2} and 
\ref{property3} (i.e., they are reflexive and ah-symmetric).

\medskip

\begin{rema}\label{rem:approx_from1_to_21}
  Theorem  \ref{thm:21rwithP3} is stronger than the direct consequence of
  Theorem \ref{thmrwithP3} that can be obtained via 
  norm inequalities. Specifically, for $x \!\in \!\mathbb{R}^m$, we have $\|x\|_2 \!\leq\! \|x\|_1 \!\leq \!\sqrt{m}\|x\|_2$ \cite[Eq. 2.2.5]{GVL1996}. Then, for  $H\!\in\!\mathbb{R}^{n\times m}$, we have $\|H\|_{2,1} \!=\! \textstyle\sum_{i = 1}^n \|H_{i\cdot}\|_2  \!\leq\!  \textstyle\sum_{i = 1}^n \|H_{i\cdot}\|_1 
\! =\! \|H\|_1$,
 and 
 $\|H\|_{1}\!= \!\textstyle\sum_{i = 1}^n \|{H}_{i\cdot}\|_1 \!\leq\!  \sqrt{m}\textstyle\sum_{i = 1}^n \|H_{i\cdot}\|_2
\! =\! \sqrt{m}\|H\|_{2,1}$\,.
Then, for $H\!\in\!\mathbb{R}^{n\times m}$ constructed as in Theorem \ref{thmrwithP3} we have 
 $\|H\|_{2,1} \!\leq\! \|H\|_1 \!\leq\! r(1+\epsilon) \|H^1_{opt}\|_1 \!\leq\! r(1+\epsilon) \|H_{opt}^{2,1}\|_{1}\!\leq\! r\sqrt{m}(1+\epsilon) \|H_{opt}^{2,1}\|_{2,1}$\,.  
That dependence on $m$ is concerning, in the  context of least-squares applications, where $m$ can be huge.
\end{rema}

\medskip

To prove Theorem \ref{thm:21rwithP3}, we will work with the convex-optimization
formulation for a 2,1-norm minimizing generalized inverse,
\begin{align}\label{prob:C}\tag{\mbox{$P_1^{2,1}$}}
&\min\left\{ \|H\|_{2,1}~:~ AHA = A\right\}, 
\end{align}
and its dual
\begin{align}\label{prob:dualC}\tag{\mbox{$D_1^{2,1}$}}
\max\left\{\tr(A^\top W)~:~\left\|(A^\top W A^\top)_{i\cdot}\right\|_2 \leq 1, ~   i=1,\dots,n\right\}.
\end{align}

Next, we present some technical results that will be used to prove  our main result.

\begin{lemm}\label{thm:feasibleE}
    Let $T$ be an ordered subset of $r$ elements from $\{1,\dots,n\}$ and $\hat{A}:= A[\cdot,T]$ be the $m\! \times \!r$ submatrix of an $m \!\times \!n$ matrix $A$ formed by columns $T$, and $\rank(\hat{A})\!=\!r$.  There exists an $r\!\times\! r$ matrix $E$ with $\| E_{i\cdot}\hat{A}^\top\|_{2} \!=\! 1$ and $E_{ii} =  \|\hat{A}^{\dagger}_{ i\cdot}\|_2\,$ for $i \in \{1,\dots,r\}$.
\end{lemm}

\begin{proof}
Suppose that $\hat{A} = U\Sigma V^\top$ is the compact singular value decomposition of $\hat{A}$, where $U\in\mathbb{R}^{m\times r}$, $\Sigma, V\in\mathbb{R}^{r\times r}$, and $\Sigma$ is diagonal.
Let $\hat{V}=\Sigma^{-1}V^\top = [z_1,\dots, z_r]$ ($\Sigma^{-1} V^\top \mathbf{e}_i = z_i$), and $E_{ij} = \frac{z_i^\top z_j}{\|z_i\|_2}$ ($E = M\hat{V}^\top\hat{V} = MV\Sigma^{-2}V^\top$, where $M$ is a diagonal matrix with $M_{ii} = \frac{1}{\|z_i\|_2}$).
We know that $\hat{A}^{\dagger} = V\Sigma^{-1} U^\top$. 
Then
$E_{ii} = \frac{z_i^\top z_i}{\|z_i\|_2} = \|z_i\|_2=\|\Sigma^{-1} V^\top \mathbf{e}_i\|_2=\|U\Sigma^{-1} V^\top \mathbf{e}_i\|_2=\|\hat{A}^{\dagger}_{ i\cdot}\|_2\,$, where the second last equality follows from orthogonality of $U$.
Also, as $\hat{A}^\top = V\Sigma U^\top$, we have $E_{i\cdot}\hat{A}^\top $ $=$ $ \frac{1}{\|z_i\|_2}\mathbf{e}_i^\top V\Sigma^{-2}V^\top V\Sigma U^\top $ $=$ $ \frac{1}{\|z_i\|_2}\mathbf{e}_i^\top V\Sigma^{-1} U^\top $ $=$ $ \frac{1}{\|z_i\|_2}\hat{A}^{\dagger}_{ i\cdot}$.
Because $\|z_i\|_2 = \|\hat{A}^{\dagger}_{ i\cdot}\|_2\,$, we have $\| E_{i\cdot}\hat{A}^\top\|_{2}=1$.
\end{proof}


\begin{lemm}\label{lem:AWE}
Let $T$ be an ordered subset of $r$ elements from $\{1,\dots,n\}$ and $\hat{A}:=A[\cdot,T]$ be the $m \times r$ submatrix of an $m \times n$ matrix $A$ formed by columns $T$, and $\rank(\hat{A})=r$. Let $E$ be an $r\times r$ matrix such that  $\| E_{i\cdot}\hat{A}^\top\|_{2} = 1$ and $E_{ii} =  \|\hat{A}^{\dagger}_{ i\cdot}\|_2$\,, for $i \in \{1,\dots,r\}$ (which exists by Lemma \ref{thm:feasibleE}). There exists an $m\times n$ matrix $W$ such that
$
\hat{A}^\top W A^\top = E\hat{A}^\top
$
 and  $
\tr(A^\top W)
=\|\hat{A}^\dagger\|_{2,1}$\,.

\end{lemm} 
\begin{proof}   
Let $S$ be an ordered subset of $r$ elements from $\{1,\dots,m\}$ such that $\tilde{A}:=\hat{A}[S,\cdot]$ is a nonsingular $r\times r$ submatrix of $\hat{A}$ formed by rows $S$.
Let $\hat{W}$ be an $r\times r$ matrix and $W$ be an $m\times n$ matrix with all elements equal to zero, except the ones in rows $S$ and columns $T$ which are given by the respective elements in $\hat{W}$. Choose $\hat{W}:=\tilde{A}^{-\top}E$, 
then we have $\hat{A}^\top W A^\top = \tilde{A}^\top \hat{W}\hat{A}^\top = E \hat{A}^\top$ and 
$
\tr(A^\top W)
=\tr(\tilde{A}^\top \hat{W}) = \tr(E) = \|\hat{A}^\dagger\|_{2,1}$\,.
\end{proof}

We now proceed to 
prove Theorem \ref{thm:21rwithP3}.

\begin{proof} 
We will construct a dual-feasible solution with objective value $\textstyle\frac{1}{r(1+\epsilon)}\| H\|_{2,1}$\,. By weak duality, we will then have $\textstyle\frac{1}{r(1+\epsilon)}\| H\|_{2,1}\le \|H_{opt}\|_{2,1}$\,.

By Lemma  \ref{lem:AWE}, we can always choose $W$ such that
$
\textstyle\hat{A}^\top W A^\top=E\hat{A}^\top
$
and 
$
\textstyle\tr(A^\top W)
= \tr(E) = \|\hat{A}^\dagger\|_{2,1}=\|H\|_{2,1}~,
$ where 
  $E$ is any given $r \times r$ matrix  such that  $\| E_{i\cdot}\hat{A}^\top\|_{2} = 1$ and $E_{ii} =  \|\hat{A}^{\dagger}_{ i\cdot}\|_2\,$, for $i \in \{1,\dots,r\}$.

So it is sufficient to demonstrate that $\|(A^\top W A^\top)_{i\cdot}\|_{2}\le r(1+\epsilon)$ for $i = 1,\dots,n$ (then $\textstyle\frac{1}{r(1+\epsilon)}W$ is dual feasible), and $\textstyle\tr{\left(A^\top\left(\textstyle\frac{1}{r(1+\epsilon)}W\right)\right)} = \textstyle\frac{1}{r(1+\epsilon)}\|H\|_{2,1}$\,.

First, it is clear that
$
\| (\hat{A}^\top W A^\top)_{i\cdot}\|_{2} = \| E_{i\cdot} \hat{A}^\top\|_{2} = 1\leq r(1+\epsilon),
$ for $i \in T$.
Now, consider any column $\hat{b}$ of $ A[\cdot,N\setminus T]$.
Because $\rank (\hat{A})=r$, we have that $\hat{b}=\hat{A}\beta$, for some $\beta\in\mathbb{R}^r$, which implies that $\tilde{b}=\tilde{A}\beta$, where 
$\tilde{b}: = \hat b[S]$.
By Cramer's rule, where $\tilde{A}_i(\tilde{b})$ is $\tilde{A}$ with column $i$ replaced by $\tilde{b}$, we have
$
|\beta_i|=\frac{|\det(\tilde{A}_i(\tilde{b}))|}{|\det(\tilde{A})|}\le 1+\epsilon,
$
because $\tilde{A}$ is a $(1+\epsilon)$-local maximizer for the absolute determinant of $A[S,\cdot]$. Therefore
\begin{align*}
\|\hat{b}^\top W A^\top\|_{2}&= \|\beta^\top \hat{A}^\top W A^\top \|_{2}
= \|\beta^\top E\hat{A}^\top\|_{2}
= \|\textstyle\sum_{i = 1}^r \beta_i \cdot E_{i\cdot} \hat{A}^\top\|_2 \\
&\le\textstyle\sum_{i=1}^r\|\beta_i\cdot E_{i\cdot} \hat{A}^\top\|_2
= \textstyle\sum_{i=1}^r|\beta_i|\cdot\|E_{i\cdot} \hat{A}^\top\|_2\le r(1+\epsilon),
\end{align*}
where the first inequality comes from the triangle inequality.
\end{proof}


\begin{rema}\label{rem:criterion}
We wish to  emphasize that even though the local search of 
Definition \ref{def:localsearch_ahsym} does not directly consider either the (vector) 1-norm or the 
2,1-norm in its operation, Theorems   \ref{thmrwithP3} and 
\ref{thm:21rwithP3}
still provide approximation guarantees for both. Furthermore, 
we showed in \cite[Example A.1(3)]{XFLPsiam}, that the approximation ratio  established in Theorem \ref{thmrwithP3} is essentially the best possible for the 1-norm. Using the same example, we can also show that the  approximation ratio established in Theorem  \ref{thm:21rwithP3}  for the 2,1-norm is also essentially best possible.
\end{rema}

\medskip

As emphasized in Remark \ref{rem:criterion}, although our goal is to minimize the 2,1-norm over ah-symmetric reflexive generalized inverses, the local-search procedure presented above uses instead, the absolute determinant of  $r\times r$ nonsingular submatrices of $A$ as a criterion for improving  the  constructed solution. The advantage of the criterion used is twofold: it leads to the approximation result shown in Theorem \ref{thm:21rwithP3} and to a very efficient implementation based on  rank-1 update of the determinant. 

Nevertheless, a natural further investigation concerns how the results compare to the solution obtained by the local search modified to use as a criterion for improving  the solution, the actual 2,1-norm of $H$. More specifically,  the criterion
for improving  the  column-block solution $H$ constructed by Theorem \ref{thm:ahconstruction}, is modified to be the decrease in its 2,1-norm, or equivalently,
the decrease in the 2,1-norm of the M-P pseudoinverse of the $m \times r$  full column-rank submatrix of $A$ being
considered. To evaluate how much the 2,1-norm of the M-P pseudoinverse of the submatrix changes when
each column of $A[\cdot, T]$ is replaced by a given column  $\gamma$  of $A[\cdot, N \!\setminus\! T ]$, we use the next result.

\begin{prop}\label{remnormcolpinv}
Let  $A:=(a_1,  \ldots, a_j\mathbin{,}\ldots\mathbin{,} a_r) \mathbin{\in} \mathbb{R}^{m\times r}$ with $\rank(A)=r$ and $\gamma\in \mathcal{R}(A)$. Let $\bar{A}$ be the matrix obtained by replacing the $j^{th}$ column of
$A$ by $\gamma$, and $v=(v_1,\ldots,v_j,\ldots,v_r)^\top:=A^{\dagger}\gamma$. If $v_j\neq 0$, let 
$\bar{v}:=\textstyle\frac{1}{v_j}(-v_1,\ldots, -v_{j-1},1, -v_{j+1},\ldots,-v_r)^\top.
$ 
 Then
\[
\textstyle\|\bar{A}^\dagger\|_{2,1} = 
|\bar{v}_j|\cdot {\|A_{j\cdot}^\dagger\|_2}+ 
\sum_{\substack{i=1\\i\neq j}}^r  
\left(\|A^{\dagger}_{i\cdot}\|_2^2  +  2\bar{v}_i A^\dagger_{i\cdot}(A^\dagger_{j\cdot})^\top + \bar{v}_i^2\|A^{\dagger}_{j\cdot}\|_2^2\right)^{1/2}\,.
\]
\end{prop}
\begin{proof}
Because $\gamma\in \mathcal{R}(A)$, we have $Av=AA^{\dagger}\gamma=\gamma$. Then, 
$\bar{A} = A\Theta$, where $\Theta:= (e_1, \ldots, e_{j-1},\  v, \ e_{j+1},\ldots, e_r)$ is non-singular, and we can verify that $\rank(\bar{A})=\rank(A)=r$ and that
\[
\bar{A}\bar{A}^{\dagger}=\bar{A}(\bar{A}^\top\bar{A})^{-1}\bar{A}^\top=A\Theta(\Theta^\top A^\top  A\Theta)^{-1}\Theta^\top A^\top =A( A^\top  A)^{-1} A^\top = AA^{\dagger}.
\]
Then, we have
\begin{align*}
&\bar{A}\bar{A}^{\dagger}=A\Theta\bar{A}^{\dagger}\Rightarrow AA^{\dagger}=A\Theta\bar{A}^{\dagger}\Rightarrow A^{\dagger}AA^{\dagger}=A^{\dagger}A\Theta\bar{A}^{\dagger} \Rightarrow A^{\dagger}=\Theta\bar{A}^{\dagger},  
\end{align*}
where the last implication holds because $A$ has full column rank.

 Then, $\bar{A}^{\dagger} = \Theta^{-1} A^{\dagger}$, with $\Theta^{-1}= (e_1, \ldots, e_{j-1},\  \bar{v}, \ e_{j+1},\ldots, e_r)$, or equivalently, 
 $
 \bar{A}^\dagger_{i\cdot} =  A^\dagger_{i\cdot} + \bar{v}_iA^\dagger_{j\cdot}\,,$ for $i\neq j$, and $\bar{A}^\dagger_{j\cdot} =  \bar{v}_jA^\dagger_{j\cdot}\,.$ So, we finally have 
\begin{align*}
 \|\bar{A}^\dagger_{i\cdot}\|_2 =\left((A^\dagger_{i\cdot} + \bar{v}_iA^\dagger_{j\cdot})(A^\dagger_{i\cdot} + \bar{v}_iA^\dagger_{j\cdot})^\top\right)^{1/2}
 =\left(\|A^\dagger_{i\cdot}\|_2^2 + 2 \bar{v}_iA^\dagger_{i\cdot}(A^\dagger_{j\cdot})^\top + \bar{v}_i^2\|A^\dagger_{j\cdot}\|_2^2\right)^{1/2},
 \end{align*}
 for $i\neq j$, and
 $\|\bar{A}^\dagger_{j\cdot}\|_2 = |\bar{v}_j|\cdot \|A^\dagger_{j\cdot}\|_2\,.$
 The result follows.
\end{proof}

\begin{rema}
    Computing $A^\dagger_{i\cdot}(A^\dagger_{j\cdot})^\top$ for a given $j$ and all $i\in \{1,\ldots,r\}$ to update the 2,1-norm of the M-P pseudoinverse, as described in Proposition \ref{remnormcolpinv}, at every iteration of our local-search procedure, is still quite  time consuming. To address this, we note that $\bar{A}^\dagger_{i\cdot} =  A^\dagger_{i\cdot} + \bar{v}_iA^\dagger_{j\cdot}\,,$ for $i\neq j$, and $\bar{A}^\dagger_{j\cdot} =  \bar{v}_jA^\dagger_{j\cdot}\,.$ Therefore, defining $W := A^\dagger (A^\dagger)^\top$ we note that it is possible to  compute  $\bar{A}^\dagger (\bar{A}^\dagger)^\top$ by  the following.
    \begin{align*}
       { \left(\bar{A}^\dagger (\bar{A}^\dagger)^\top\right)_{i\ell} = \begin{cases}
        W_{i\ell} +  \bar{v}_i\bar{v}_\ell\|{A}^\dagger_{j\cdot}\|_{2}^2 +  W_{ij}\bar{v}_\ell + W_{j\ell}\bar{v}_i\,,\quad& i\neq j,\, \ell \neq j;\\
        \bar{v}_j\left(W_{i\ell} + \|{A}^\dagger_{j\cdot}\|_{2}^2\bar{v}_i\right),&i\neq j,\, \ell = j;\\
        \bar{v}_j^2\|{A}^\dagger_{j\cdot}\|_{2}^2,&i=j,\, \ell = j.
    \end{cases}
    }
    \end{align*}
    The  update above drastically improves upon computing 
   the  $A^\dagger_{i\cdot}(A^\dagger_{j\cdot})^\top$ from scratch.    
\end{rema}



\section{Targeting row sparsity with ADMM}\label{sec:admm20}

It is a common approach to induce row-sparsity for a matrix by minimizing its 2,1-norm, in order to take advantage of the convexity of the problem addressed (see \cite{yuan2006model,meier2008group,obozinski2008high,liu2012multi} and \cite[Sec. 6.4.2]{boyd2011distributed}, for example). However, we can observe from our numerical experiments, that the ah-symmetric reflexive generalized inverse of a given matrix $A\in\mathbb{R}^{m\times n}$,  obtained with this approach, has in general significantly more nonzero rows than the solution of the local-search procedure described in \S\ref{sec:LS}, for which the number of nonzero rows, or equivalently, the 2,0-norm,  is $r:=\mbox{rank}(A)$, the least possible number. On the other side, we observe that the solution obtained by the local search has in general much larger 2,1-norm than a 2,1-norm minimizing    
generalized inverse. Aiming at matrices  with both nice features of having small 2,1- and 2,0-norms, we now seek  an ah-symmetric reflexive generalized inverse  with 2,1-norm smaller  than the local-search solution  and 2,0-norm smaller  than the  2,1-norm minimizing solution. Inspired by ideas in  \cite{LiuYu2023},
we will present  ADMM algorithms to compute an ah-symmetric reflexive generalized inverses for $A$, with the 2,0-norm limited to $\gamma:=\omega r+(1-\omega)\|H_{opt}^{2,1}\|_{2,0}$\,, 
where $0 < \omega <1$, and $H_{opt}^{2,1}$ an optimal solution of \ref{prob:barmin21norm}. 
More specifically, we will present ADMM algorithms for  \emph{nonconvex}  problems of obtaining a solution $\bar Z$ of $\|V_1D^{-1} + V_2 Z\|_{2,0}\leq \gamma$. 
Initially, aiming at a more efficient application of the ADMM method,  we do not consider the minimization of the 2,1-norm, so the problem addressed is a nonconvex feasibility problem. Then, we investigate the impact on the solutions  by minimizing the 2,1-norm subject to the same nonconvex feasible set.
Our goal is to compare through numerical experiments, the solutions obtained by the  ADMM algorithms proposed next, with each other and also to the local-search solutions and the solutions of \ref{prob:barmin21norm}. 

We wish to note that
the ADMM approach that  we 
develop, for handling 
an inequality-constrained 
minimization problem in one variable 
by introducing a second variable,  
an indicator function on the 
second variable, and 
linear linking constraints
is a known scheme (at a high level);
see~\cite[introductory passages of Section~5, and Sections 5.2, 6.2 and 9.1]{boyd2011distributed}, for example.

\subsection{ADMM for limited 2,0-norm}\label{subsec:ADMM20}
By introducing a variable $E \in \mathbb{R}^{n \times r}$, we  reformulate the feasibility problem of obtaining a solution $\bar Z$ of $\|V_1D^{-1} + V_2 Z\|_{2,0}\leq \gamma$ as
\begin{equation}\label{prob:20converteda}
\min\{ \mathcal{I}_{\mathcal{M}}(E) ~:~ E =  V_1D^{-1} + V_2 Z\},
\end{equation}
where  the  indicator function $\mathcal{I}_{\mathcal{M}}(\cdot)$ is defined by
\begin{align*}
    \mathcal{I}_{\mathcal{M}}(X) := \begin{cases}
    0, & X\in\mathcal{M};\\ 
    +\infty, &X\notin \mathcal{M}.
\end{cases}
\end{align*}
for the set $\mathcal{M}:=\{X\in\mathbb{R}^{n\times r}~:~\|X\|_{2,0} \leq \gamma \}$. 

The augmented Lagrangian function associated to \eqref{prob:20converteda} is
\begin{align*}
\mathcal{L}_\rho(Z,E,\Lambda)\!:=\!\textstyle
\mathcal{I}_{\mathcal{M}}(E) \!+\! \frac{\rho}{2} \left\|V_1D^{-1} \!+\! V_2 Z \!-\! E\!+\! \Lambda\right\|^2_F - 
  \frac{\rho}{2}\left\|\Lambda\right\|_F^2\,,
\end{align*}
where $\rho >0$ is the penalty parameter and $\Lambda \in \mathbb{R}^{n \times r}$ is the scaled Lagrangian multiplier. We will apply the ADMM method to \eqref{prob:20converteda}, by iteratively solving, for $k=0,1,\ldots,$
\begin{align}
    &Z^{k+1}:=\textstyle\argmin_Z ~ \mathcal{L}_\rho(Z,E^{k},\Lambda^k),\label{eq:Zmin20normsubpa}\\
    &E^{k+1}:=\textstyle\argmin_E ~ \mathcal{L}_\rho(Z^{k+1},E,\Lambda^k),\label{eq:Emin20normsubpa}\\
    &\textstyle\Lambda^{k+1}:=\Lambda^{k} + V_1D^{-1} + V_2 Z^{k+1} - E^{k+1}.\nonumber
\end{align}
Subproblem \eqref{eq:Zmin20normsubpa} is exactly the same as \eqref{eq:Zmin21normsubpa}. Next, we detail how to solve 
\eqref{eq:Emin20normsubpa}.

\medskip

\noindent{\bf{Update of \texorpdfstring{$E$}{E}:}}
\begin{align*}
    E^{k+1}
    :=\textstyle\argmin_{E}\{\mathcal{I}_{\mathcal{M}}(E) \!+\! \frac{\rho}{2} \left\|E-Y\right\|^2_F\},
\end{align*}
where $Y:=V_1D^{-1} + V_2 Z^{k+1} + \Lambda^k$.
 Although $\mathcal{M}$ is a nonconvex set, the exact solution for the above subproblem can be efficiently computed (see \cite[Section 9.1]{boyd2011distributed}).   It is given by the projection of $Y$ onto $\mathcal{M}$, which we will denote in the following by $\Pi_\mathcal{M}(Y)$. $\Pi_{\mathcal{M}}(Y)$ keeps the  rows of $Y$ with the $\gamma$ largest 2-norms and zeros out the other rows. 


\medskip

\noindent{\bf{Initialization:}} 
To initialize the variables, we follow the approach described in \S\ref{subsec:admm1} and  \S\ref{subsec:ADMM21}. But  because of the nonconvexity of \eqref{prob:20converteda}, instead of considering its dual problem to initialize $\Lambda$, we simply set $\Lambda^0:=0$. As in \S\ref{subsec:ADMM21}, we set $E^0 := V_1D^{-1} + \Lambda^0$.

\medskip

\noindent{\bf{Stopping criteria:}}
Considering the  nonconvexity of  our feasibility problem, formulated as \eqref{prob:20converteda}, we stop the ADMM  when a feasible solution is found, i.e., when  $\|V_1D^{-1}\! + V_2 Z^{k}\|_{2,0}\!\leq\! \gamma$, (or when the time limit is reached). 

\medskip

\noindent{\bf{Pseudocode:}}
In Algorithm \ref{alg:admm20norm}, we present the ADMM algorithm for \eqref{prob:20converteda}\,.

\begin{algorithm2e}[!ht]
\footnotesize{
\KwIn{ $A\in \mathbb{R}^{m\times n}$,
$\Lambda^0 \in \mathbb{R}^{n\times r}$, $E^0 \in \mathbb{R}^{n\times r}$, 
  $0<\omega<1$, $r:=\mbox{rank}(A)$, $\|H_{opt}^{2,1}\|_{2,0}$\,.}
\KwOut{$H\in \mathbb{R}^{n\times m}$.}
$U,\Sigma,V := \texttt{svd}(A)$, $k:=0$\;
 Get $U_1,V_1,V_2,D^{-1}$ from $U,\Sigma,V$\;
$\mathcal{M}:=\{X\in\mathbb{R}^{n\times r}~:~\|X\|_{2,0} \leq \omega r +(1-\omega)\|H_{opt}^{2,1}\|_{2,0} \}$\;
\While {not converged} 
{
$J := E^k -V_1D^{-1} - \Lambda^k$\;
$Z^{k+1} := V_2^\top J$\;
$Y := V_1D^{-1} + V_2 Z^{k+1} + \Lambda^k$\;
$E^{k+1}:=  \Pi_{\mathcal{M}}(Y)$\;
$\Lambda^{k+1} := \Lambda^k + V_1D^{-1} + V_2 Z^{k+1} - E^{k+1}$\;
$k:=k+1$\;
}
$H:=V_1D^{-1}U_1^\top + V_2Z^kU_1^\top$\;
\caption{ADMM for \eqref{prob:20converteda}\label{alg:admm20norm} (ADMM$_{2,0}$)}
\hypertarget{algadmm20norm}{}
}
\end{algorithm2e}


\subsection{ADMM for 2,1-norm minimization subject to limited 2,0-norm}\label{subsec:admm2120}

Next, we introduce  variable $E \in \mathbb{R}^{n \times r}$ in    \ref{prob:barmin21norm} and restrict it  to solutions with 2,0-norm limited to $\gamma$\,, 
 reformulating it as
\begin{equation}\label{prob:admm_mat2120}
\min\limits_{E,Z}\left\{ \|E\|_{2,1}~:~
 E =  V_1D^{-1} + V_2 Z,~ \|E\|_{2,0}\leq \gamma\right\}.
\end{equation}
Then, we reformulate  \eqref{prob:admm_mat2120} as
\begin{equation}\label{prob:2120converteda}
\min\limits_{E,Z}\left\{\|E\|_{2,1} +  \mathcal{I}_{\mathcal{M}}(E) ~:~ E =  V_1D^{-1} + V_2 Z\right\}.
\end{equation}
The augmented Lagrangian function associated to \eqref{prob:2120converteda} is

\begin{align*}
\mathcal{L}_\rho(Z,E,\Lambda)\!:=\!\textstyle
 \|E\|_{2,1} +  \mathcal{I}_{\mathcal{M}}(E) + \frac{\rho}{2}\left\|V_1D^{-1} + V_2Z -E + \Lambda\right\|^2_F  - 
  \frac{\rho}{2}\left\|\Lambda\right\|_F^2, 
\end{align*}
where $\rho >0$ is the penalty parameter and $\Lambda \in \mathbb{R}^{n \times r}$  is the scaled Lagrangian multiplier. We will apply the ADMM method to \eqref{prob:2120converteda}, by iteratively solving, for $k=0,1,\ldots$,
\begin{align}
    &Z^{k+1}:=\textstyle\argmin_Z ~ \mathcal{L}_\rho(Z,E^{k},\Lambda^k),\label{eq:Zmin2120normsubpa}\\
    &E^{k+1}:=\textstyle\argmin_E ~ \mathcal{L}_\rho(Z^{k+1},E,\Lambda^k),\label{eq:Emin2120normsubpa}\\
    &\textstyle\Lambda^{k+1}:=\Lambda^{k} + V_1D^{-1} + V_2 Z^{k+1} - E^{k+1}.\nonumber
\end{align}
Subproblem \eqref{eq:Zmin2120normsubpa} is exactly the same as \eqref{eq:Zmin21normsubpa}. Next, we detail how to solve \eqref{eq:Emin2120normsubpa}.

\medskip

\noindent{\bf{Update of 
\texorpdfstring{$E$}{B}:}} 
To update $E$,  we consider  subproblem \eqref{eq:Emin2120normsubpa}, that is

\begin{align}
    E^{k+1}&:=\textstyle\argmin_{E \in \mathbb{R}^{n \times r}}\left\{ \|E\|_{2,1}
    +\mathcal{I}_{\mathcal{M}}(E) + \textstyle\frac{\rho}{2}\left\|V_1D^{-1} + V_2 Z^{k+1} -E + \Lambda^k\right\|^2_F\right\}\nonumber\\
    &=\textstyle \argmin_{E \in \mathcal{M}} \left\{\|E\|_{2,1}
    + \textstyle\frac{\rho}{2}\left\|E-Y\right\|^2_F\right\},\label{subprobB}
\end{align}
where $Y := V_1D^{-1} + V_2 Z^{k+1} + \Lambda^k$.

\begin{theo}\label{lem:2120_projection}
Let $\tilde{E}^{k+1} \in \mathbb{R}^{n \times r}$ be the optimal solution to the unconstrained version of  subproblem  \eqref{subprobB}, that is, 
\begin{equation}\label{unconst}
\tilde{E}^{k+1}:=\argmin\{\|E\|_{2,1}
    + \textstyle\frac{\rho}{2}\left\|E-Y\right\|^2_F\}.
    \end{equation}
    Let 
    \begin{equation}\label{eq:g(i)2120}
        g(i) := \textstyle\frac{\rho}{2}\|Y_{i\cdot}\|_2^2 - \|\tilde{E}^{k+1}_{i\cdot}\|_2 - \textstyle\frac{\rho}{2}\|\tilde{E}^{k+1}_{i\cdot} - Y_{i\cdot}\|_2^2 \,,
    \end{equation}
    for $i\in N:=\{1,\dots,n\}$. Let $\phi$ be the permutation of the indices in $N$ such that $g(\phi_1) \geq g(\phi_2) \geq \dots \geq g(\phi_n)$. 
    Then, an optimal solution $E^{k+1} \in \mathbb{R}^{n \times r}$ of \eqref{subprobB} is given by \begin{equation}\label{updateE_2120_lemma_zero_out}
    E^{k+1}_{\phi_i\cdot} := \begin{cases}
    \tilde{E}^{k+1}_{\phi_i\cdot} ,\quad &\text{if } i \leq \gamma;
    \\
    0,&\text{otherwise.}
\end{cases}
\end{equation}
\end{theo}

\begin{proof}
From \eqref{updateE}, we see that 
\begin{itemize}
    \item[($i$)]\label{case1}   if $1\!/\!\rho <  \|Y_{i\cdot}\|_2$  , then $\tilde{E}^{k+1}_{i\cdot}\neq 0$, and we can verify that $g(i)=\textstyle\frac{\rho}{2}\left(\|Y_{i\cdot}\| - \frac{1}{\rho}\right)^2>0$,
    \item[($ii$)]\label{case2} if $1\!/\!\rho \geq  \|Y_{i\cdot}\|_2$  , then $\tilde{E}^{k+1}_{i\cdot}=0$, and we can verify that  $g(i)=0$.
    \end{itemize}
\smallskip
From ($i$) and ($ii$), we see that  ${E}^{k+1}\!=\!\tilde{E}^{k+1}$  if $\|\tilde{E}^{k+1}\|_{2,0} \!\leq\! \gamma$, and the statement of the theorem trivially follows. 
Therefore, in the following, we assume that $\|\tilde{E}^{k+1}\|_{2,0}\! >\! \gamma$. Note that in this case we have  that  $\|\tilde{E}^{k+1}_{\phi_i\cdot}\|_2 > 0$, for all $i\leq \gamma$, so  
\begin{equation}\label{normE}
\|{E}^{k+1}_{\phi_i\cdot}\|_2 > 0, \mbox{ for all } i\leq \gamma \mbox{ and }  \|E^{k+1}\|_{2,0} = \gamma.
\end{equation}

Next, we will demonstrate
that there is no feasible solution to \eqref{subprobB} with a better objective value than that of $E^{k+1}$.
Let us  suppose, by contradiction, that there exists  a solution $X \in \mathbb{R}^{n \times r}$ with $\|X\|_{2,0} \!\leq\! \gamma$,  such that
    \begin{align}
        & \textstyle\|X\|_{2,1} + \textstyle\frac{\rho}{2}\|X - Y\|_F^2 < \|E^{k+1}\|_{2,1} + \textstyle\frac{\rho}{2}\|E^{k+1} - Y\|_F^2\Leftrightarrow \nonumber \\
        &\qquad \textstyle\sum_{i \in N } \|X_{i\cdot}\|_2 + \textstyle\frac{\rho}{2}\|X_{i\cdot} - Y_{i\cdot}\|_2^2 < \sum_{i \in N}\|E^{k+1}_{i\cdot}\|_2 + \textstyle\frac{\rho}{2}\|E^{k+1}_{i\cdot} - Y_{i\cdot}\|_2^2\,.\label{sumcontrad}
    \end{align}
    In this case, there exists  $\hat\jmath \in N$, such that     \begin{equation}\label{contrad}
     \textstyle\|X_{\hat\jmath\cdot}\|_2+ \textstyle\frac{\rho}{2}\|X_{\hat\jmath\cdot} - Y_{\hat\jmath\cdot}\|_2^2  < \|E^{k+1}_{\hat\jmath\cdot}\|_2 + \textstyle\frac{\rho}{2}\|E^{k+1}_{\hat\jmath\cdot} - Y_{\hat\jmath\cdot}\|_2^2\,.
     \end{equation} 
Because  $\tilde{E}^{k+1}$ is optimal  and $X$ is feasible for  \eqref{unconst},
we must have 
        \begin{equation}\label{bestsol}
           \|X_{\hat\jmath\cdot}\|_2 + \textstyle\frac{\rho}{2}\|X_{\hat\jmath\cdot} - Y_{\hat\jmath\cdot}\|_2^2\geq \textstyle\|\tilde{E}^{k+1}_{\hat\jmath\cdot}\|_2 + \textstyle\frac{\rho}{2}\|\tilde{E}^{k+1}_{\hat\jmath\cdot} - Y_{\hat\jmath\cdot}\|_2^2\,,
        \end{equation}
         otherwise  we would obtain a better solution than $\tilde{E}^{k+1}$ for \eqref{unconst}, by replacing its $\hat\jmath$-th row with $X_{\hat\jmath\cdot}$\,.
         
          From \eqref{contrad} and \eqref{bestsol}, 
        we can see that   $E^{k+1}_{\hat\jmath\cdot} \!\neq\! 
    \tilde{E}^{k+1}_{\hat\jmath\cdot}$. Thus, from \eqref{updateE_2120_lemma_zero_out}, we have $E^{k+1}_{\hat\jmath\cdot}\!=\!0$, 
    and  \eqref{contrad} reduces~to 
        \[
        \textstyle \|X_{\hat\jmath\cdot}\|_2+ \textstyle\frac{\rho}{2}\|X_{\hat\jmath\cdot} - Y_{\hat\jmath\cdot}\|_2^2 <
         \textstyle\frac{\rho}{2}\|Y_{\hat\jmath\cdot}\|_2^2\,.
         \]
        Note that we cannot have  $X_{\hat\jmath\cdot} = 0$ in the inequality above, hence  we have $\|X_{\hat\jmath\cdot}\|_2 > 0$  and $\|E^{k+1}_{\hat\jmath\cdot}\|_2=0$. 
        
        Now, we recall that  $\|E^{k+1}\|_{2,0} = \gamma$ 
        and  $\|X\|_{2,0} \leq \gamma$, therefore for each  $\hat\jmath$ that satisfies \eqref{contrad},  there must exist a distinct $\hat\ell \in N \setminus \{\hat\jmath\}$  such that  
        $\|X_{\hat\ell\cdot}\|_2 = 0$ and $\|E^{k+1}_{\hat\ell\cdot}\|_2 > 0$. 
Moreover, because 
        $\|E^{k+1}_{\hat\ell\cdot}\|_{2} > 0$ and $\|E^{k+1}_{\hat\jmath\cdot}\|_{2} =0$, we see from \eqref{updateE_2120_lemma_zero_out} and \eqref{normE} that $\hat\ell\in\{\phi_1,\phi_2, \dots, \phi_\gamma\}$  and $\hat\jmath\in\{\phi_{\gamma+1},\phi_{\gamma+2}, \dots, \phi_n\}$, so 
        $g(\hat\ell) \geq g(\hat\jmath)$,  that is
  \begin{equation}\label{eq:lemma_2120_eq_EkEj2}
             \textstyle\frac{\rho}{2}\|Y_{\hat\ell\cdot}\|_2^2 - \|\tilde{E}^{k+1}_{\hat\ell\cdot}\|_2 - \frac{\rho}{2}\|\tilde{E}^{k+1}_{\hat\ell\cdot} - Y_{\hat\ell\cdot}\|_2^2 \geq  
            \frac{\rho}{2}\|Y_{\hat\jmath\cdot}\|_2^2 - \|\tilde{E}^{k+1}_{\hat\jmath\cdot}\|_2 - \frac{\rho}{2}\|\tilde{E}^{k+1}_{\hat\jmath\cdot} - Y_{\hat\jmath\cdot}\|_2^2\,.
\end{equation}    
        Also,  from \eqref{bestsol}, we can see that 
\begin{equation}\label{contrad3}
        \textstyle\frac{\rho}{2}\|Y_{\hat\jmath\cdot}\|_2^2 - \|\tilde{E}^{k+1}_{\hat\jmath\cdot}\|_2 - \frac{\rho}{2}\|\tilde{E}^{k+1}_{\hat\jmath\cdot} - Y_{\hat\jmath\cdot}\|_2^2 \geq 
        \textstyle\frac{\rho}{2}\|Y_{\hat\jmath\cdot}\|_2^2 - \|X_{\hat\jmath\cdot}\|_2 - \frac{\rho}{2}\|X_{\hat\jmath\cdot} - Y_{\hat\jmath\cdot}\|_2^2\,.
\end{equation}
        From \eqref{eq:lemma_2120_eq_EkEj2} and \eqref{contrad3}, we have that
\[
         \textstyle \frac{\rho}{2}\|Y_{\hat\ell\cdot}\|_2^2 - \|\tilde{E}^{k+1}_{\hat\ell\cdot}\|_2 - \frac{\rho}{2}\|\tilde{E}^{k+1}_{\hat\ell\cdot} - Y_{\hat\ell\cdot}\|_2^2 \geq \textstyle\frac{\rho}{2}\|Y_{\hat\jmath\cdot}\|_2^2 - \|X_{\hat\jmath\cdot}\|_2 - \frac{\rho}{2}\|X_{\hat\jmath\cdot} - Y_{\hat\jmath\cdot}\|_2^2\,.
\]
         We recall that    $X_{\hat\ell\cdot} = 0$ and $E^{k+1}_{\hat\jmath\cdot} = 0$. Moreover,  $\|E^{k+1}_{\hat\ell\cdot}\|_2 > 0$, thus $E_{\hat\ell\cdot}^{k+1} = \tilde{E}^{k+1}_{\hat\ell\cdot}$. Then,  the last inequality  is equivalent to
         \begin{align*}
        & \|X_{\hat\ell\cdot}\|_2 +\textstyle\frac{\rho}{2}\|X_{\hat\ell\cdot} \!-\!Y_{\hat\ell\cdot}\|_2^2 - \|{E}^{k+1}_{\hat\ell\cdot}\|_2 -\textstyle\frac{\rho}{2}\|{E}^{k+1}_{\hat\ell\cdot} \!-\! Y_{\hat\ell\cdot}\|_2^2 \geq  \|{E}^{k+1}_{\hat\jmath\cdot} \|_2 +\textstyle\frac{\rho}{2}\|{E}^{k+1}_{\hat\jmath\cdot} \!-\!Y_{\hat\jmath\cdot}\|_2^2 - \|X_{\hat\jmath\cdot}\|_2 -\textstyle\frac{\rho}{2}\|X_{\hat\jmath\cdot} \!-\! Y_{\hat\jmath\cdot}\|_2^2\,,
         \end{align*}
         which we rewrite as 
         \begin{equation}\label{ineqjell}
         \textstyle\sum_{i \in \{\hat\jmath,\hat\ell\}} \|{E}^{k+1}_{i\cdot} \|_2 + \frac{\rho}{2}\|{E}^{k+1}_{i\cdot} -Y_{i\cdot}\|_2^2 - \|X_{i\cdot}\|_2 -\frac{\rho}{2}\|X_{i\cdot} -Y_{i\cdot}\|_2^2\leq 0.
         \end{equation}

         Taking into account the above, we see that for each $\hat\jmath$ that satisfies \eqref{contrad}, there is a distinct $\hat\ell:=\hat\ell(\hat\jmath)$, such that \eqref{ineqjell} holds for $(\hat\jmath, \hat{\ell}(\hat\jmath))$. Let $\hat{N}$ be the set of all indices $\hat\jmath\in N$ that satisfy \eqref{contrad} and the corresponding indices $\hat\ell(\hat\jmath)$. Then we have  
  \begin{align*}
        &\textstyle\sum_{i \in N}\|E^{k+1}_{i\cdot}\|_2 + \textstyle\frac{\rho}{2}\|E^{k+1}_{i\cdot} - Y_{i\cdot}\|_2^2 - \textstyle\sum_{i \in N }\left( \|X_{i\cdot}\|_2 + \textstyle\frac{\rho}{2}\|X_{i\cdot} - Y_{i\cdot}\|_2^2\right) \\
        &\quad =  \textstyle\sum_{i \in \hat{N}} \|{E}^{k+1}_{i\cdot} \|_2 + \frac{\rho}{2}\|{E}^{k+1}_{i\cdot} -Y_{i\cdot}\|_2^2 - \|X_{i\cdot}\|_2 -\frac{\rho}{2}\|X_{i\cdot} -Y_{i\cdot}\|_2^2 \\
        &\quad \quad \quad + \textstyle\sum_{i \in N\setminus\hat{N}} \|{E}^{k+1}_{i\cdot} \|_2 + \frac{\rho}{2}\|{E}^{k+1}_{i\cdot} -Y_{i\cdot}\|_2^2 - \|X_{i\cdot}\|_2 -\frac{\rho}{2}\|X_{i\cdot} -Y_{i\cdot}\|_2^2 \\
        &\quad \leq 0 \,,
    \end{align*}
    which contradicts \eqref{sumcontrad}, showing that there is no better solution than $E^{k+1}$ to \eqref{subprobB}.
\end{proof}

\begin{rema}
    It is possible to verify that the result of Theorem \ref{lem:2120_projection}
     still holds if subproblem \eqref{subprobB} is replaced by the more general problem 
    \[
    E^{k+1}:=\textstyle \argmin_{E \in \mathcal{M}}\sum_{i=1}^{n} f(E_{i\cdot})
    + \frac{\rho}{2}\left\|E-Y\right\|^2_F,
\]
where $f:\mathbb{R}^r\rightarrow \mathbb{R}$ is nonnegative  with $f(\delta)=0$ if and only if $\delta=0$. 
For example, we could have the 2,1-norm in the objective function of subproblem \eqref{subprobB} replaced by the 1-norm or the square of the Frobenius norm.
\end{rema}

\medskip

\begin{coro}\label{updE2120}
    Let $\tau$ be the permutation of the
indices in $N:=\{1,\ldots,n\}$, such that $\|Y_{\tau_1\cdot}\|_2 \geq  \|Y_{\tau_2\cdot}\|_2 \geq \cdots \geq \|Y_{\tau_n\cdot}\|_2$\,. Then an optimal solution of \eqref{subprobB} is given by 
\begin{equation}\label{updateE_2120}
    E^{k+1}_{\tau_i\cdot} := \begin{cases}
    \frac{ \|Y_{\tau_i\cdot}\|_2 - 1\!/\!\rho}{\|Y_{\tau_i\cdot}\|_2} Y_{\tau_i\cdot}\, ,\quad &\text{if } 1\!/\!\rho < \|Y_{\tau_i\cdot}\|_2 \text{ and } i \leq \gamma \,;
    \\
    0,&\text{otherwise.}
\end{cases}
\end{equation}
\end{coro}
\begin{proof}
Let us consider $\tilde{E}^{k+1}$ as the optimal solution of the unconstrained version of \eqref{subprobB} 
defined in Proposition \ref{lem:closedformula21norm}, that is  
\begin{equation}
    \tilde E^{k+1}_{\tau_i\cdot} := \begin{cases}
    \frac{ \|Y_{\tau_i\cdot}\|_2 - 1\!/\!\rho}{\|Y_{\tau_i\cdot}\|_2} Y_{\tau_i\cdot}\, ,\quad &\text{if } 1\!/\!\rho < \|Y_{\tau_i\cdot}\|_2  \,;
    \\
    0,&\text{otherwise.}
\end{cases}
\end{equation}

Considering the notation and statement of Theorem \ref{lem:2120_projection}, it suffices to prove that $g(\tau_1) \geq g(\tau_2) \geq \dots \geq g(\tau_n)$.
 For all $i$, such that  $1\!/\!\rho < \|Y_{\tau_i\cdot}\|_2$\,, we have
   \begin{align*}
       g(\tau_i) &=\textstyle\frac{\rho}{2}\|Y_{\tau_i\cdot}\|_2^2 - \left\|\left(\frac{ \|Y_{\tau_i\cdot}\|_2 - 1\!/\!\rho}{\|Y_{\tau_i\cdot}\|_2}\right) Y_{\tau_i\cdot}\right\|_2 - \frac{\rho}{2}\left\| \left(\frac{ \|Y_{\tau_i\cdot}\|_2 - 1\!/\!\rho}{\|Y_{\tau_i\cdot}\|_2}\right) Y_{\tau_i\cdot} - Y_{\tau_i\cdot}\right\|_2^2\\
       &=\textstyle\frac{\rho}{2}\|Y_{\tau_i\cdot}\|_2^2 - \left(\frac{ \|Y_{\tau_i\cdot}\|_2 - 1\!/\!\rho}{\|Y_{\tau_i\cdot}\|_2}\right) \left\|Y_{\tau_i\cdot}\right\|_2 - \frac{\rho}{2}\left(1 - \frac{ \|Y_{\tau_i\cdot}\|_2 - 1\!/\!\rho}{\|Y_{\tau_i\cdot}\|_2}\right)^2\left\|  Y_{\tau_i\cdot}\right\|_2^2\\
       &=\textstyle\frac{\rho}{2}\|Y_{\tau_i\cdot}\|_2^2 - \|Y_{\tau_i\cdot}\|_2 
 + \textstyle\frac{1}{\rho} - \frac{1}{2\rho}
 =\textstyle\left(\sqrt{\frac{\rho}{2}} \|Y_{\tau_i\cdot} \|_2 - \sqrt{\frac{1}{2\rho}}\;\right)^2.
   \end{align*}
We note that, because $1\!/\!\rho < \|Y_{\tau_i\cdot}\|_2$\,, we have $\textstyle\sqrt{\frac{\rho}{2}} \|Y_{\tau_i\cdot} \|_2 - \sqrt{\frac{1}{2\rho}} > 0$. Then, for any pair of indices $j,\ell$ with $j<\ell$, $1\!/\!\rho < \|Y_{\tau_j\cdot}\|_2$\,, and $1\!/\!\rho < \|Y_{\tau_\ell\cdot}\|_2$\,,  as we have   $\|Y_{\tau_j\cdot}\|_2 \geq  \|Y_{\tau_\ell\cdot}\|_2\,$, we also have  $g(\tau_j) \geq g(\tau_{\ell})>0$. Moreover, for all $i$, such that $1\!/\!\rho \geq \|Y_{\tau_i\cdot}\|_2$\,, we have $g(\tau_i)=0$, confirming that $g(\tau_1) \geq g(\tau_2) \geq \dots \geq g(\tau_n)$. Finally, from the ordering  $\|Y_{\tau_1\cdot}\|_2 \geq  \|Y_{\tau_2\cdot}\|_2 \geq \dots \geq \|Y_{\tau_n\cdot}\|_2$\,, it is clear  that $\|\tilde{E}^{k+1}_{\tau_1\cdot}\|_2 \geq \|\tilde{E}^{k+1}_{\tau_2\cdot}\|_2\geq \cdots\geq \|\tilde{E}^{k+1}_{\tau_n\cdot}\|_2$\,, completing the proof.
\end{proof}

\medskip

\noindent{\bf{Initialization of the variables:}} 
We initialize the variables with the solution obtained by the ADMM for 2,1-norm minimization described in \S\ref{subsec:ADMM21}.

\medskip

\noindent{\bf{Stopping criteria:}}
We adopt the  stopping criterion described in \S\ref{subsec:ADMM21}, additionally  requiring that  $\|H\|_{2,0}\leq \gamma$.

\medskip

\noindent{\bf{Pseudocode:}}
In Algorithm \ref{alg:admm2120norm}, we present the ADMM algorithm for \eqref{prob:2120converteda}\,.

\begin{algorithm2e}[!ht]
\footnotesize{
\KwIn{ $A\in \mathbb{R}^{m\times n}$,
$\Lambda^0 \in \mathbb{R}^{n\times r}$, $E^0 \in \mathbb{R}^{n\times r}$, 
  $\rho > 0$, $0<\omega<1$, $r:=\mbox{rank}(A)$, $\|H_{opt}^{2,1}\|_{2,0}$\,.}
\KwOut{$H\in \mathbb{R}^{n\times m}$.}
$U,\Sigma,V := \texttt{svd}(A)$, $k:=0$\;
 Get $U_1,V_1,V_2,D^{-1}$ from $U,\Sigma,V$\;
$\mathcal{M}:=\{X\in\mathbb{R}^{n\times r}~:~\|X\|_{2,0} \leq \gamma :=\omega r +(1-\omega)\|H_{opt}^{2,1}\|_{2,0} \}$\;
\While {not converged} 
{
$J := E^k -V_1D^{-1} - \Lambda^k$\;
$Z^{k+1} := V_2^\top J$\;
$Y := V_1D^{-1} + V_2 Z^{k+1} + \Lambda^k$\;
$\tau:=$  permutation of
indices in $\{1,\dots,n\}$, such that $\|Y_{\tau_1\cdot}\|_2 \geq  \|Y_{\tau_2\cdot}\|_2 
\geq \!\cdots\! \geq \|Y_{\tau_n\cdot}\|_2$ \;
\For{$i = 1,\dots,n$}{
    \bf{if} $\|Y_{\tau_i\cdot}\|_2> 1\!/\!\rho~ \& ~i\leq \gamma$ \bf{then} 
   $E^{k+1}_{\tau_i\cdot} := \frac{\|Y_{\tau_i\cdot}\|_2 - 1\!/\!\rho}{\|Y_{\tau_i\cdot}\|_2} Y_{\tau_i\cdot}$; ~   \mbox{(see Cor. \ref{updE2120})}{\color{white}\;}
    \bf{else}
    $E^{k+1}_{\tau_i\cdot} := 0$\;
}
$\Lambda^{k+1} := \Lambda^k + V_1D^{-1} + V_2 Z^{k+1} - E^{k+1}$\;
$k:=k+1$\;
}
$H:=V_1D^{-1}U_1^\top + V_2Z^kU_1^\top$\;
\caption{ADMM for \eqref{prob:2120converteda} (ADMM$_{2,1/0}$) \label{alg:admm2120norm}}
\hypertarget{algadmm2120norm}{}
}
\end{algorithm2e}


\section{Numerical experiments}\label{sec:Exper}

We constructed test instances 
of varying sizes for our numerical experiments using the MATLAB function \texttt{sprand} to  randomly generate $m\!\times\! n$  matrices $A$ with  rank $r$, as described in \cite[Section 2.1]{FLPXjogo}. 
The instances considered in \cite{FLPXjogo} were too small for our experiments, but we considered the results of this previous work by selecting the formulations for our norm minimization problems as those in which \texttt{Gurobi} and \texttt{MOSEK} performed best, namely \ref{prob:barmin1norm} and \ref{prob:barmin21norm}.
We divided our instances into two categories related to $m$ with $n: = 0.5m$, $r := 0.25m$; small instances with $m := 100,\, 200,\dots,\,500$ (S1, S2,$\dots$, S5) and large instances with $m := 1000,\, 2000,\dots,\, 5000$ (L1, L2,$\dots$, L5).    

We ran our experiments on `zebratoo', a
32-core machine (running Windows Server 2022 Standard):
two Intel Xeon Gold 6444Y processors running at 3.60GHz, with 16 cores each, and 128 GB of memory. We consider $10^{-5}$ as the tolerance to distinguish nonzero elements. The symbol `*'  in column `Time' of our tables, indicates   that the problem was not solved to optimality because the time limit   was reached, and the symbol `\skull' indicates that we ran out of  memory. We set a time limit of 2 hours for solving each instance. 
We coded our  algorithms   in \texttt{Julia} v.1.10.0.


In our first numerical experiment, we 
compare the solutions of \hyperlink{algadmm1norm}{ADMM$_{1}$}  for  \ref{prob:barmin1norm} with the solutions of \texttt{Gurobi}, and  the solutions of \hyperlink{algadmm21norm}{ADMM$_{2,1}$}  for  \ref{prob:barmin21norm} with the solutions of \texttt{MOSEK}. 
We used \texttt{Gurobi} v.11.0.0 to solve  \ref{prob:barmin1norm} as a linear-optimization problem:

\begin{equation*}\label{gur_LP}
\min_{\newatop{F\in\mathbb{R}^{n\times m},}{Z\in\mathbb{R}^{(n-r)\times r}}}\!\!\left\{ 
\textstyle\sum_{i=1}^{n}\sum_{j=1}^m F_{ij}:~
F -  V_2ZU_1^\top\geq V_1 D^{-1}U_1^\top\,,~F +  V_2ZU_1^\top\geq -V_1 D^{-1}U_1^\top\,\right\},
\end{equation*}
and we used \texttt{MOSEK}  v.10.1.21 to solve  \ref{prob:barmin21norm}
a second-order-cone optimization problem:

\begin{equation*}\label{Mosek_socp}
\min_{\newatop{t \in {\mathbb{R}^{n}},}
{Z \in \mathbb{R}^{(n-r) \times r}}}\left\{\textstyle\sum_{i = 1}^n t_i~:~\left\|\mathbf{e}_i^\top V \begin{bmatrix}
            D^{-1}\\
            Z
        \end{bmatrix}\right\|_2\leq t_i\,, \quad i = 1,\dots,n\right\}.
\end{equation*}

Our first goal is to compare  solutions  and running times of the  ADMM algorithms to   state-of-the-art 
general-purpose convex-optimization solvers. For this comparison, in addition to running the ADMM algorithms with the stopping criteria described in \S\ref{subsec:admm1} and \S\ref{subsec:ADMM21}, where the residual tolerances depend not only on the parameters $\epsilon^{\mbox{\scriptsize abs}}$ and $\epsilon^{\mbox{\scriptsize rel}}$, but also on the values of the variables that are dynamically updated throughout the iterations (see (\ref{sc1r}--\ref{sc1s}), (\ref{sc2r}--\ref{sc2s})), 
we also run the ADMM algorithms  with a fixed tolerance $\epsilon$ for the Frobenius norms of the primal and dual residuals. With this fixed tolerance, we have a more fair comparison with the solvers that also work with a fixed tolerance.  In the following, we refer to this second version of ADMM$_1$ and ADMM$_{2,1}$ with a fixed tolerance $\epsilon$, respectively, by ADMM$_1^\epsilon$ and  ADMM$_{2,1}^\epsilon$\,. 

Our second goal is  to compare the solutions and running times of ADMM$_1$ (and ADMM$_1^\epsilon$)   to ADMM$_{2,1}$ (and ADMM$_{2,1}^\epsilon$). With these comparisons, we verify whether the different norms used to induce sparsity and row sparsity are effective and how much we actually lose in sparsity and gain in row sparsity  when applying ADMM$_{2,1}$ compared to ADMM$_1$\,. 

In Table \ref{tab:solvers}, we show from the first to the last column, the instance, the method adopted, 
the different norms 
of the solutions and the (elapsed) running time (in seconds). In the second column, the labels identifying the methods  are: \texttt{Gurobi}, \hyperlink{algadmm1norm}{ADMM$_{1}$}\,, and  ADMM$_1^\epsilon$   applied to solve \ref{prob:barmin1norm}\,; \texttt{MOSEK},  \hyperlink{algadmm21norm}{ADMM$_{2,1}$}\,, and ADMM$_{2,1}^\epsilon$  applied to solve \ref{prob:barmin21norm}. 

The following parameters were used. 
\begin{itemize}
    \item \texttt{Gurobi}: optimality and feasibility tolerances of $10^{-4}$; 
     \item \hyperlink{algadmm1norm}{ADMM$_{1}$}\,: $\epsilon^{\mbox{\scriptsize abs}} = \epsilon^{\mbox{\scriptsize rel}} := 10^{-4}$, $\rho := 3$;
    \item ADMM$_1^\epsilon$\,: $\epsilon := 10^{-4}$, $\rho := 3$;
 \item \texttt{MOSEK}: optimality and feasibility tolerances of $10^{-5}$; 
   \item \hyperlink{algadmm21norm}{ADMM$_{2,1}$}\,: $\epsilon^{\mbox{\scriptsize abs}}=\epsilon^{\mbox{\scriptsize rel}}:=10^{-7}$, $\rho := 1$.
    \item ADMM$_{2,1}^\epsilon$\,: $\epsilon := 10^{-5}$, $\rho := 1$;
\end{itemize}

We note that we initially tried to select the values of $\epsilon^{\text{abs}}$, $\epsilon^{\text{rel}}$ and $\epsilon$,  all equal to $10^{-4}$. However, when we ran our experiments with \texttt{MOSEK} using tolerance $10^{-4}$, we observed solutions with much higher 2,0-norms than the solutions of  ADMM$_{2,1}^\epsilon$ (probably because \texttt{MOSEK} is using interior-point methods). To get a fair comparison between \texttt{MOSEK} and ADMM$_{2,1}^\epsilon$ with  comparable solutions in terms of row-sparsity, we selected the smaller tolerance of $10^{-5}$ for both. Finally, the small tolerance of $10^{-7}$ was selected for \hyperlink{algadmm21norm}{ADMM$_{2,1}$} to highlight the fast convergence of this algorithm even for small tolerances, which was not the case for \hyperlink{algadmm1norm}{ADMM$_{1}$}\,.

\begin{table}[!ht]
\centering
\scriptsize
\tabcolsep=3pt
\caption{Comparison between ADMM and solvers}\label{tab:solvers}
\begin{tabular}{l|l|rrrrr}
Inst.               & Method                & $\|H\|_1$ & $\|H\|_0$ & $\|H\|_{2,1}$ & $\|H\|_{2,0}$ & Time (sec) \\   \hline \strut
\multirow{6}{*}{S1} & \texttt{Gurobi}       & 194.28    & 3734      & 39.14         & 44            &6.06 \\
                    & ADMM$_{1}^\epsilon$   & 194.28    & 3733      & 39.14         & 44            & 1.67       \\
                    & \hyperlink{algadmm1norm}{ADMM$_{1}$}            & 194.29    & 4152      & 39.13         & 50            & 0.52       \\
                    & \texttt{MOSEK}        & 211.60    & 3824      & 36.10         & 42            & 1.43       \\
                    & ADMM$_{2,1}^\epsilon$ & 211.60    & 3819      & 36.10         & 39            & 0.08       \\
                    &\hyperlink{algadmm21norm}{ADMM$_{2,1}$}         & 211.60    & 3819      & 36.10         & 39            & 0.08       \\   \hline  \strut
\multirow{6}{*}{S2} & \texttt{Gurobi}       & 539.70    & 14943     & 79.91         & 90            & 119.21 \\ 
                    & ADMM$_{1}^\epsilon$   & 539.71    & 14927     & 79.91         & 90            & 71.58      \\
                    & \hyperlink{algadmm1norm}{ADMM$_{1}$}            & 539.79    & 16804     & 79.89         & 99            & 5.30       \\
                    & \texttt{MOSEK}        & 595.24    & 16222     & 73.41         & 84            & 0.38       \\
                    & ADMM$_{2,1}^\epsilon$ & 595.24    & 16219     & 73.41         & 83            & 0.01       \\
                    &\hyperlink{algadmm21norm}{ADMM$_{2,1}$}         & 595.24    & 16219     & 73.41         & 83            & 0.02       \\ \hline   \strut
\multirow{6}{*}{S3} & \texttt{Gurobi}       & 868.19    & 33756     & 119.40        & 133           &2295.54 \\ 
                    & ADMM$_{1}^\epsilon$   & 868.19    & 33702     & 119.40        & 134           & 309.21     \\
                    & \hyperlink{algadmm1norm}{ADMM$_{1}$}            & 868.39    & 38665     & 119.40        & 150           & 12.67      \\
                    & \texttt{MOSEK}        & 972.46    & 34723     & 109.98        & 118           & 29.96      \\
                    & ADMM$_{2,1}^\epsilon$ & 972.47    & 34729     & 109.98        & 117           & 0.24       \\
                    &\hyperlink{algadmm21norm}{ADMM$_{2,1}$}         & 972.47    & 34729     & 109.98        & 117           & 0.33       \\ \hline   \strut
\multirow{6}{*}{S4} & \texttt{Gurobi}       & -         & -         & -             & -             & *          \\
                    & ADMM$_{1}^\epsilon$   & 1339.90   & 64174     & 163.73        & 194           & 870.73     \\
                    & \hyperlink{algadmm1norm}{ADMM$_{1}$}            & 1340.28   & 73374     & 163.73        & 200           & 19.94      \\
                    & \texttt{MOSEK}        & 1504.61   & 67338     & 150.07        & 193           & 2.33       \\
                    & ADMM$_{2,1}^\epsilon$ & 1504.64   & 66956     & 150.07        & 171           & 0.26       \\
                    &\hyperlink{algadmm21norm}{ADMM$_{2,1}$}         & 1504.64   & 66956     & 150.07        & 171           & 0.33       \\ \hline   \strut
\multirow{6}{*}{S5} & \texttt{Gurobi}       & -         & -         & -             & -             & *          \\
                    & ADMM$_{1}^\epsilon$   & 1870.64   & 100758    & 208.31        & 245           & 1030.52    \\
                    & \hyperlink{algadmm1norm}{ADMM$_{1}$}            & 1871.23   & 114769    & 208.30        & 250           & 27.39      \\
                    & \texttt{MOSEK}        & 2102.71   & 106021    & 190.28        & 237           & 5.48       \\
                    & ADMM$_{2,1}^\epsilon$ & 2102.77   & 105312    & 190.28        & 217           & 0.31       \\
                    &\hyperlink{algadmm21norm}{ADMM$_{2,1}$}         & 2102.77   & 105312    & 190.28        & 217           & 0.48       \\ \hline   \strut
\multirow{6}{*}{L1} & \texttt{Gurobi}       & -         & -         & -             & -             & *          \\
                    & ADMM$_{1}^\epsilon$   & 4886.76   & 409821    & 425.60        & 492           & 3178.79    \\
                    & \hyperlink{algadmm1norm}{ADMM$_{1}$}            & 4888.92   & 465117    & 425.68        & 500           & 50.17      \\
                    & \texttt{MOSEK}        & 5637.29   & 437653    & 384.94        & 480           & 40.90      \\
                    & ADMM$_{2,1}^\epsilon$ & 5637.67   & 435820    & 384.94        & 444           & 1.03       \\
                    &\hyperlink{algadmm21norm}{ADMM$_{2,1}$}         & 5637.67   & 435820    & 384.94        & 444           & 0.58       \\ \hline   \strut
\multirow{6}{*}{L2} & \texttt{Gurobi}       & -         & -         & -             & -             & *          \\
                    & ADMM$_{1}^\epsilon$   & 11781.79  & 1622951   & 853.10        & 987           & *          \\
                    & \hyperlink{algadmm1norm}{ADMM$_{1}$}            & 11789.25  & 1848689   & 853.23        & 1000          & 252.32     \\
                    & \texttt{MOSEK}        & 14098.24  & 1757540   & 765.78        & 978           & 443.30     \\
                    & ADMM$_{2,1}^\epsilon$ & 14099.58  & 1746738   & 765.78        & 891           & 2.14       \\
                    &\hyperlink{algadmm21norm}{ADMM$_{2,1}$}         & 14099.58  & 1746741   & 765.78        & 891           & 2.26       \\ \hline   \strut
\multirow{6}{*}{L3} & \texttt{Gurobi}       & -         & -         & -             & -             & *          \\
                    & ADMM$_{1}^\epsilon$   & 20459.62  & 3672370   & 1306.86       & 1488          & *          \\
                    & \hyperlink{algadmm1norm}{ADMM$_{1}$}            & 20474.95  & 4183302   & 1307.02       & 1500          & 563.40     \\
                    & \texttt{MOSEK}        & -         & -         & -             & -             & \skull     \\
                    & ADMM$_{2,1}^\epsilon$ & 24904.76  & 4014677   & 1160.70       & 1367          & 4.25       \\
                    &\hyperlink{algadmm21norm}{ADMM$_{2,1}$}         & 24904.76  & 4014676   & 1160.70       & 1367          & 4.05       \\ \hline   \strut
\multirow{6}{*}{L4} & \texttt{Gurobi}       & -         & -         & -             & -             & *          \\
                    & ADMM$_{1}^\epsilon$   & 30029.31  & 6474511   & 1738.97       & 1995          & *          \\
                    & \hyperlink{algadmm1norm}{ADMM$_{1}$}            & 30054.56  & 7407976   & 1738.99       & 2000          & 946.15     \\
                    & \texttt{MOSEK}        & -         & -         & -             & -             & \skull     \\
                    & ADMM$_{2,1}^\epsilon$ & 36984.26  & 7075074   & 1547.08       & 1813          & 7.65       \\
                    &\hyperlink{algadmm21norm}{ADMM$_{2,1}$}         & 36984.26  & 7075080   & 1547.08       & 1813          & 7.23       \\ \hline   \strut
\multirow{6}{*}{L5} & \texttt{Gurobi}       & -         & -         & -             & -             & *          \\
                    & ADMM$_{1}^\epsilon$   & 40716.32  & 10285553  & 2191.30       & 2498          & *          \\
                    & \hyperlink{algadmm1norm}{ADMM$_{1}$}            & 40751.46  & 11584101  & 2191.83       & 2498          & 1506.57    \\
                    & \texttt{MOSEK}        & -         & -         & -             & -             & \skull     \\
                    & ADMM$_{2,1}^\epsilon$ & 50516.23  & 11223830  & 1947.08       & 2297          & 14.40      \\
                    &\hyperlink{algadmm21norm}{ADMM$_{2,1}$}         & 50516.22  & 11223836  & 1947.08       & 2297          & 12.91     
\end{tabular}
\end{table}


From the results in Table \ref{tab:solvers}, we see that with the application of \texttt{Gurobi} to our linear-optimization model for \ref{prob:barmin1norm}\,, we can only solve  the three smallest instances.  ADMM$_{1}^{\epsilon}$ converges for all small instances and for one large instance as well. When using the looser dynamic stopping criterion, we have convergence of \hyperlink{algadmm1norm}{ADMM$_{1}$} for all instances tested. ADMM$_{1}^{\epsilon}$ and \texttt{Gurobi}  converge to  very similar solutions when both converge, but ADMM$_{1}^{\epsilon}$ converges much faster, showing the advantage of the application of the ADMM algorithm over our application of \texttt{Gurobi} to construct  1-norm minimizing ah-symmetric reflexive generalized inverses.   \hyperlink{algadmm1norm}{ADMM$_{1}$}\,, with the dynamic stopping criterion, converges much faster than  ADMM$_{1}^{\epsilon}$
for all tested instances and, comparing the solutions of both algorithms, we see only an increase of at most 0.1\% in the 1-norm of the solutions obtained by
 \hyperlink{algadmm1norm}{ADMM$_{1}$}\,, 
 although the 0-norms increase by about 10\%.

Analyzing now the statistics for ADMM$_{2,1}^{\epsilon}$ and ADMM$_{2,1}$\,, we see that both algorithms  converge for all instances tested to solutions with the same 2,1-norms, and both take less than  15 seconds on the largest instance.  The ADMM for the 2,1-norm is robust and quickly converges to high-precision solutions.
\texttt{MOSEK} did not solve the three largest instances due to lack of memory, and for the other instances, \texttt{MOSEK} took much longer than the ADMM algorithms and obtained solutions with the same 2,1-norms.
The results in Table \ref{tab:solvers} show  high superiority of the ADMM algorithm proposed  to construct 2,1-norm minimizing ah-symmetric reflexive generalized inverses, when compared to our application of \texttt{MOSEK}.  
We note  from the 2,0-norms of the solutions,  that the minimization of the 2,1-norm is in fact a more effective strategy to induce row sparsity than the minimization of the 1-norm. Moreover,  we always obtain fewer nonzero rows with \hyperlink{algadmm21norm}{ADMM$_{2,1}$} than with \texttt{MOSEK}. 

\begin{figure}[!ht]%
    \centering
    \subfloat[]{{\includegraphics[width=0.48\textwidth]{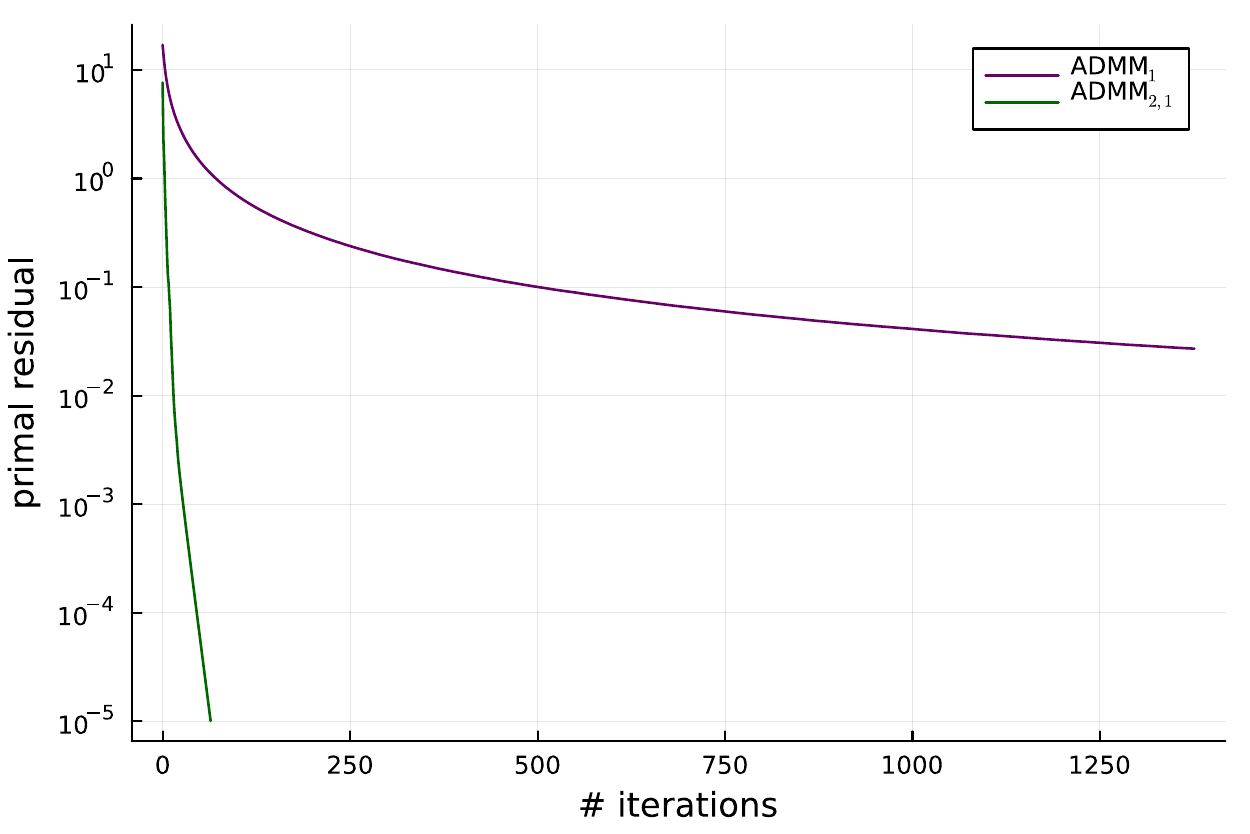} }}%
    ~
    \subfloat[]{{\includegraphics[width=0.48\textwidth]{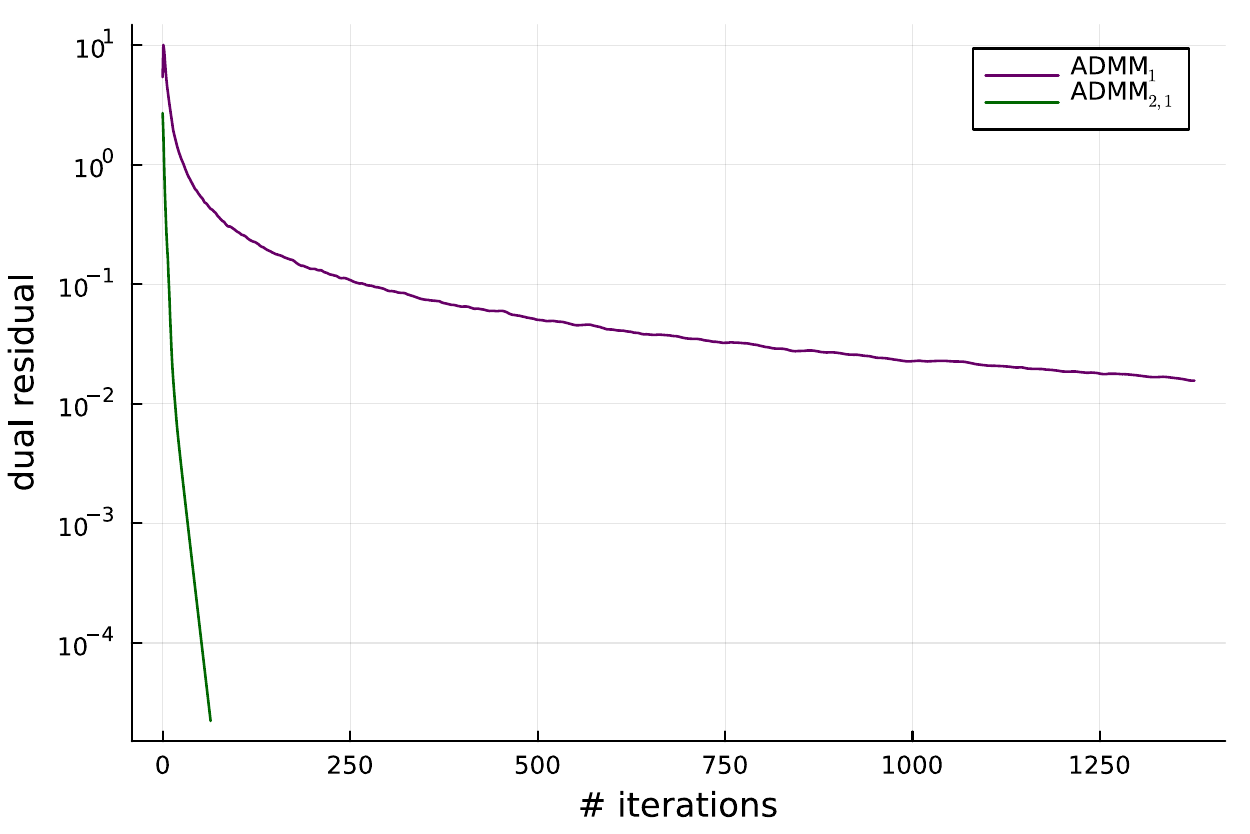} }}%
    \hfill
     \caption{Convergence: \protect\hyperlink{algadmm1norm}{ADMM$_{1}$}  vs.\! \protect\hyperlink{algadmm21norm}{ADMM$_{2,1}$} ; $m,n,r=1000,500,250$ (Instance L1) } 
    \label{fig:compare_bounds1}%
\end{figure}

Finally, we can observe from the statistics in Table \ref{tab:solvers}, that ADMM$_{2,1}$ converges much faster than \hyperlink{algadmm1norm}{ADMM$_{1}$}\,. In Fig. \ref{fig:compare_bounds1}, we illustrate the convergence of these algorithms, showing that the number of iterations  for \hyperlink{algadmm21norm}{ADMM$_{2,1}$} is also much lower. 
Although Fig. \ref{fig:compare_bounds1} corresponds to instance L1 only, it illustrates the typical behavior in our experiments.

We also note from Table \ref{tab:solvers} that, as expected,  ADMM$_1^\epsilon$ obtains the  sparsest solutions among the four ADMM algorithms tested, but at a high computational cost. Therefore, if our goal is to compute solutions quickly,  we should note that although  \hyperlink{algadmm21norm}{ADMM$_{2,1}$} aims at row sparsity, it is actually more effective than \hyperlink{algadmm1norm}{ADMM$_{1}$} at obtaining sparse solutions.

Finally, we note that our goal when setting the parameters of the solvers was to achieve the best possible performance. In the case of \texttt{Gurobi}, this was achieved by letting \texttt{Gurobi} choose the algorithm. We also experimented with selecting specific algorithms, but \texttt{Gurobi} did not perform well. From the \texttt{Gurobi} log for an instance, we can see that it starts with a barrier method and then uses a dual simplex method, but this can change depending on the difficulty of the instance for \texttt{Gurobi}. 

For \texttt{MOSEK}, we can see that it recognizes each instance as a conic optimization problem, for which \texttt{MOSEK} employs a primal-dual interior point algorithm. From \texttt{MOSEK}'s log for some instances, we can see that it performs many fewer iterations than \hyperlink{algadmm21norm}{ADMM$_{2,1}$} and ADMM$_{2,1}^\epsilon$\,, but each iteration of \texttt{MOSEK} takes much longer than the iterations of the ADMM algorithm. We can then conclude that, for the instances considered, solving the linear system of equations to obtain the Newton direction in the primal-dual interior point method iterations is much more computationally expensive than updating the variables in \hyperlink{algadmm21norm}{ADMM$_{2,1}$} with the closed-form formulas presented in Algorithm \ref{alg:admm21norm}.

In our second numerical experiment, our goal is to verify the efficiency of the ADMM   algorithms described in \S\ref{sec:admm20}  
(\hyperlink{algadmm20norm}{ADMM$_{2,0}$} 
and  \hyperlink{algadmm2120norm}{ADMM$_{2,1/0}$}), in computing ah-symmetric reflexive generalized inverses within the target  2,0-norm given by $\omega r+(1-\omega)\|H_{opt}^{2,1}\|_{2,0}$\,, where $H_{opt}^{2,1}$ is the solution obtained by \hyperlink{algadmm21norm}{ADMM$_{2,1}$}\,. We also  aim at observing the trade-off between the norms of the ah-symmetric reflexive generalized inverses obtained by these algorithms while varying $\omega$ in the interval $(0,1)$.  We note that for $\omega=0$, the target 2,0-norm is the 2,0-norm of the 2,1-norm minimizing generalized inverse $H_{opt}^{2,1}$\,, while for $\omega=1$, the target 2,0-norm is $r$, the smallest possible  2,0-norm of an ah-symmetric reflexive generalized inverse. Then, increasing $\omega$ in the interval  $(0,1)$, we intend to construct ah-symmetric reflexive generalized inverses with decreasing 2,0-norms, and we would like to see the impact on the  1- and 2,1-norms of the matrices. We are particularly interested in verifying if the 1- and 2,1-norms of the solutions obtained by \hyperlink{algadmm20norm}{ADMM$_{2,0}$}
are too large, and how effective \hyperlink{algadmm2120norm}{ADMM$_{2,1/0}$} is in decreasing the 2,1-norm  with respect to 
\hyperlink{algadmm20norm}{ADMM$_{2,0}$}\,. We recall that the problems addressed by both algorithms are nonconvex, therefore there is no guarantee of convergence of these ADMM algorithms. The first is a feasibility problem and seeks any solution within the target 2,0-norm, while the second seeks solutions within the target 2,0-norm and minimum 2,1-norm.  
We run both ADMM algorithms for the large instances, for $\omega=0.25,0.50,0.75,0.80,0.90,0.95$, and compare the solutions obtained with  $H_{opt}^{2,1}$ and with  the local-search solutions  described in \S\ref{sec:LS}.  The first one, denoted in the following by LS, uses the absolute determinant as a criterion to improve a given solution. The linearly-independent rows and columns of $A$ used to construct its initial solution are obtained from the QR factorization of $A$. The second one, denoted in the following by LS$_{2,1}$\,, starts with the solution of LS, and try to improve it using the 2,1-norm as the criterion for improvement. In this case, the result in Proposition \ref{remnormcolpinv} is used in the implementation. 
We recall that both local-search procedures  construct  ah-symmetric reflexive generalized inverses with 2,0-norms equal to $r$.


\begin{table}[!ht]
 \scriptsize
\centering
\tabcolsep=3pt
\caption{\protect\hyperlink{algadmm20norm}{ADMM$_{2,0}$}~/ \protect\hyperlink{algadmm2120norm}{ADMM$_{2,1/0}$} for target 2,0-norm,  for various $\omega$}\label{tab:ADMM2120}
\begin{tabular}{l|l|rr|rr|rr|r|rr}
Inst.               & \multicolumn{1}{c|}{Method} & \multicolumn{2}{c|}{$\|H\|_1$ }& \multicolumn{2}{c|}{$\|H\|_0$} & \multicolumn{2}{c|}{$\|H\|_{2,1}$} & $\|H\|_{2,0}$ & \multicolumn{2}{c}{Time (sec)} \\ \hline  \strut
\multirow{9}{*}{L1} & \hyperlink{algadmm21norm}{ADMM$_{2,1}$}  & \multicolumn {2}{c|}{\phantom{1}5637.7} &   \multicolumn {2}{c|}{435820}    & \multicolumn {2}{c|}{384.9} &    444 &   \multicolumn {2}{c}{0.6} \\ 
 & $\omega\!=\!0.25$                & 6044.9 & 5796.4 & 388481 & 388547 & 393.5 & 387.1 &  395 & 0.8 & 26.9 \\ 
 & $\omega\!=\!0.50$                & 6836.2 & 6325.5 & 341730 & 341640 & 411.6 & 398.7 &  347 & 1.4 & 20.2 \\ 
 & $\omega\!=\!0.75$                & 8956.9 & 7686.8 & 293605 & 293557 & 483.2 & 441.9 &  298 & 3.2 & 26.3 \\ 
 & $\omega\!=\!0.80$                & 9759.6 & 8691.6 & 283825 & 283825 & 512.1 & 474.1 & 288 & 3.8 & 36.9 \\ 
 & $\omega\!=\!0.90$                & 12467.8 & 11164 & 265108 & 265133 & 613.0 & 566.2 & 269 & 8.1 & 42.6 \\ 
 & $\omega\!=\!0.95$                & 15223.1 & 13143.9 & 255298 & 255287 & 729.5 & 643.8 &  259 & 19.5 & 33.1 \\ 
 & LS                               & \multicolumn {2}{c|}{10797.9}    & \multicolumn {2}{c|}{246392}   & \multicolumn {2}{c|}{589.9}   & 250 &  \multicolumn {2}{c}{4.4} \\ 
 & LS$_{2,1}$                       & \multicolumn {2}{c|}{10452.9}   & \multicolumn {2}{c|}{246373}   & \multicolumn {2}{c|}{570.5}   & 250 &   \multicolumn {2}{c}{9.4} \\ \hline  \strut
\multirow{9}{*}{L2} & \hyperlink{algadmm21norm}{ADMM$_{2,1}$}  & \multicolumn {2}{c|}{14099.6} &   \multicolumn {2}{c|}{1746741} &   \multicolumn {2}{c|}{\phantom{1}765.8} &   891 &   \multicolumn {2}{c}{\phantom{1}2.3} \\ 
 & $\omega\!=\!0.25$                 & 15084.8 & 14576.0 & 1557984 & 1558091 & 781.5 & 769.2 &  793 & 2.8 & 22.0 \\ 
 & $\omega\!=\!0.50$                & 17169.9 & 16068.1 & 1367070 & 1366850 & 807.8 & 790.0 & 695 & 5.4 & 85.6 \\ 
 & $\omega\!=\!0.75$                & 23083.7 & 20739.7 & 1173280 & 1175007 & 932.4 & 883.6 & 597 & 11.6 & 93.1 \\ 
 & $\omega\!=\!0.80$                & 26747.9 & 22149.5 & 1133992 & 1137706 & 1020.3 & 915.0  & 578 & 12.6 & 83.7 \\ 
 & $\omega\!=\!0.90$                & 34401.6 & 30997.8 & 1061345 & 1061284 & 1231.7 & 1138.1  & 539 & 81.8 & 155.5 \\ 
 & $\omega\!=\!0.95$                &  51969.9 & 40647.4 & 1022152 & 1022053 & 1734.5 & 1412.2 & 519 & 273.3 & 272.9 \\ 
 & LS                               & \multicolumn {2}{c|}{32796.0}   & \multicolumn {2}{c|}{\phantom{1}984545}   &\multicolumn {2}{c|}{ 1272.5}   & 500 &   \multicolumn {2}{c}{15.2} \\ 
 & LS$_{2,1}$                       & \multicolumn {2}{c|}{31583.0}   & \multicolumn {2}{c|}{\phantom{1}984526}   & \multicolumn {2}{c|}{1232.4}   & 500 &   \multicolumn {2}{c}{66.2} \\ \hline  \strut
\multirow{9}{*}{L3} & \hyperlink{algadmm21norm}{ADMM$_{2,1}$}  & \multicolumn {2}{c|}{24904.8} &   \multicolumn {2}{c|}{4014676} &   \multicolumn {2}{c|}{1160.7} &   1367 &   \multicolumn {2}{c}{4.0} \\ 
 & $\omega\!=\!0.25$                & 26922.6 & 25896.1 & 3568117 & 3567608 & 1184.4 & 1166.2  & 1212 & 6.9 & 114.5 \\ 
 & $\omega\!=\!0.50$                & 32515.0 & 29371.3 & 3116491 & 3115482 & 1249.9 & 1206.6  & 1058 & 15.1 & 121.1 \\ 
 & $\omega\!=\!0.75$                & 45489.9 & 40631.5 & 2660981 & 2663715 & 1482.5 & 1395.7  & 904 & 28.1 & 183.0 \\ 
 & $\omega\!=\!0.80$                & 51695.7 & 45042.0 & 2568254 & 2572501 & 1609.5 & 1485.9  & 873 & 34.8 & 287.2 \\ 
 & $\omega\!=\!0.90$                & 70276.9 & 61176.6 & 2387213 & 2389950 & 2027.8 & 1826.9 & 811 & 69.1 & 283.6 \\ 
 & $\omega\!=\!0.95$                &  94905.5 & 83185.4 & 2298796 & 2298668 & 2580.9 & 2322.2 & 780 & 376.9 & 589.0 \\ 
 & LS                               & \multicolumn {2}{c|}{69338.6}   & \multicolumn {2}{c|}{2209609}   & \multicolumn {2}{c|}{2131.2}   & 750 &   \multicolumn {2}{c}{\phantom{3}26.8} \\ 
 & LS$_{2,1}$                       & \multicolumn {2}{c|}{62062.7}   & \multicolumn {2}{c|}{2208832}   & \multicolumn {2}{c|}{1974.1}   & 750 &   \multicolumn {2}{c}{380.4} \\ \hline  \strut
\multirow{9}{*}{L4} & \hyperlink{algadmm21norm}{ADMM$_{2,1}$}  & \multicolumn {2}{c|}{\phantom{1}36984.3} &   \multicolumn {2}{c|}{7075080} &   \multicolumn {2}{c|}{1547.1} &    1813 &   \multicolumn {2}{c}{7.2} \\ 
 & $\omega\!=\!0.25$                & 40413.6 & 38780.2 & 6302160 & 6301531 & 1582.1 & 1556.5  & 1609 & 9.2 & 147.5 \\ 
 & $\omega\!=\!0.50$                & 49000.8 & 44611.0 & 5519866 & 5517665 & 1665.2 & 1612.3 & 1406 & 19.6 & 225.8 \\ 
 & $\omega\!=\!0.75$                & 69070.6 & 60934.9 & 4722393 & 4727485 & 1967.3 & 1839.5 & 1203 & 35.5 & 437.8 \\ 
 & $\omega\!=\!0.80$                & 77766.4 & 68502.2 & 4559968 & 4567287 & 2117.3 & 1964.7  & 1162 & 42.6 & 423.8 \\ 
 & $\omega\!=\!0.90$                & 112003.5 & 93549.8 & 4247329 & 4250758 & 2747.0 & 2417.1  & 1081 & 103.0 & 727.3 \\ 
 & $\omega\!=\!0.95$                & 154855.3 & 129358.9 & 4086725 & 4090408 & 3592.4 & 3125.0  & 1040 & 183.5 & 864.6 \\ 
 & LS                               & \multicolumn {2}{c|}{124026.2}   & \multicolumn {2}{c|}{3932651}   & \multicolumn {2}{c|}{3162.3}   & 1000 &   \multicolumn {2}{c}{\phantom{13}93.9} \\ 
 & LS$_{2,1}$                       & \multicolumn {2}{c|}{111118.3}   & \multicolumn {2}{c|}{3932550}   & \multicolumn {2}{c|}{2902.5}   & 1000 &   \multicolumn {2}{c}{1335.4} \\ \hline  \strut
\multirow{9}{*}{L5} & \hyperlink{algadmm21norm}{ADMM$_{2,1}$}  & \multicolumn {2}{c|}{\phantom{1}50516.2} &   \multicolumn {2}{c|}{11223836} &   \multicolumn {2}{c|}{1947.1} &    2297 &   \multicolumn {2}{c}{12.9} \\   
 & $\omega\!=\!0.25$                & 55039.3 & 53092.0 & 9971973 & 9981248 & 1987.4 & 1958.6& 2035 & 10.4 & 126.2 \\ 
 & $\omega\!=\!0.50$                &  66928.5 & 61609.2 & 8713776 & 8715794 & 2091.4 & 2033.1  & 1773 & 24.3 & 681.0 \\ 
 & $\omega\!=\!0.75$                &  95752.6 & 84276.1 & 7426600 & 7439856 & 2478.9 & 2329.6  & 1511 & 48.7 & 456.6 \\ 
 & $\omega\!=\!0.80$                &  106328.0 & 93241.4 & 7172048 & 7185451 & 2636.9 & 2461.0  & 1459 & 61.0 & 668.3 \\ 
 & $\omega\!=\!0.90$                &  161475.1 & 131572.0 & 6661864 & 6672091 & 3563.9 & 3073.8 & 1354 & 136.8 & 935.5 \\ 
 & $\omega\!=\!0.95$                &   222013.5 & 185469.1 & 6413109 & 6412843 & 4653.3 & 4009.2  & 1302 & 332.1 & 1214.1 \\ 
 & LS                               & \multicolumn {2}{c|}{195014.1}   & \multicolumn {2}{c|}{\phantom{1}6160665}  & \multicolumn {2}{c|}{4401.7}   & 1250 &   \multicolumn {2}{c}{\phantom{5}192.4} \\ 
 & LS$_{2,1}$                       & \multicolumn {2}{c|}{168359.4}   & \multicolumn {2}{c|}{\phantom{1}6159949}   & \multicolumn {2}{c|}{3893.4}   & 1250 &    \multicolumn {2}{c}{5256.7} 
\end{tabular}
\end{table}

In Table \ref{tab:ADMM2120}, we have the same information presented in Table \ref{tab:solvers}. In the second column, we identify the methods addressed. For each value of $\omega$, we have in each  column, first the information for \hyperlink{algadmm20norm}{ADMM$_{2,0}$}\,, and then for \hyperlink{algadmm2120norm}{ADMM$_{2,1/0}$}\,. Both algorithms always obtain solutions with the target 2,0-norm, and therefore, there is only one result for $\|H\|_{2,0}$ for each $\omega$.
Next, we present some important details about these algorithms. 
\begin{itemize}
    \item For  \hyperlink{algadmm20norm}{ADMM$_{2,0}$} and \hyperlink{algadmm2120norm}{ADMM$_{2,1/0}$}\,,  
     we consider  $H_{i\cdot}$ to be a row of all zeros if $\|H_{i\cdot}\|_2 \leq 10^{-5}$. 
     
    \item  The  parameter $\rho$ is not used in the calculations of \hyperlink{algadmm20norm}{ADMM$_{2,0}$}\,, 
    so any $\rho >0$ could be used to define the augmented Lagrangian function considered in \S\ref{subsec:ADMM20}.
    \item For \hyperlink{algadmm2120norm}{ADMM$_{2,1/0}$}\,, we consider $\rho := 10^4$ until we obtain $\|H\|_{2,0} \leq \gamma$ and the Frobenius norm of the primal residual is less than $10^{-4}$. Then we update $\rho$ at each iteration if the Frobenius norm of the primal residual remains less than $10^{-4}$,
    by setting $\rho:=\rho/\alpha$, where $\alpha$ ranges from 2 to 1, starting at the largest value.
    The convergence of ADMM for a nonconvex problem is not guaranteed; and the selection of $\rho$ described above was important to obtain convergence of \hyperlink{algadmm2120norm}{ADMM$_{2,1/0}$} for all tested instances (see, for example, \cite[Sec. 3.4.1]{boyd2011distributed}, for a discussion of the update of $\rho$ in the ADMM algorithm). The large value of $\rho$ used to start the algorithm is related to the observed difficulty in satisfying the nonconvex constraint. Large values for $\rho$ leads to a fast convergence of the primal residual to zero,  so we were able to satisfy  $\|H\|_{2,0} \leq \gamma$. However,  with large $\rho$ we have  very slow convergence of the Frobenius norm of the dual residual to zero. Therefore,  when we have $\|H\|_{2,0} \leq \gamma$ with Frobenius norm of the primal residual less than  $10^{-4}$, we  set $\rho:=\rho/\alpha$, if the Frobenius norm of the primal residual remains less than  $10^{-4}$.  Otherwise, we iteratively decrease $\alpha$ and then update $\rho$, until the Frobenius norm of the primal residual remains less than $10^{-4}$ with the current $\rho$ or $\alpha=1$.  
     \item  For \hyperlink{algadmm2120norm}{ADMM$_{2,1/0}$}\,, in addition to using the solution of \hyperlink{algadmm21norm}{ADMM$_{2,1}$} to initialize the algorithm, we  experimented with initializing it with the solution of \hyperlink{algadmm20norm}{ADMM$_{2,0}$}\,, however, in this case the algorithm failed to converge for some instances.
    
    \item For \hyperlink{algadmm2120norm}{ADMM$_{2,1/0}$}\,, we set the parameters $\epsilon^{\mbox{\scriptsize abs}}=\epsilon^{\mbox{\scriptsize rel}}:=10^{-4}$. 
\end{itemize}

\smallskip

\begin{figure}[!ht]%
    \centering
    \subfloat[]{{\includegraphics[width=0.48\textwidth]{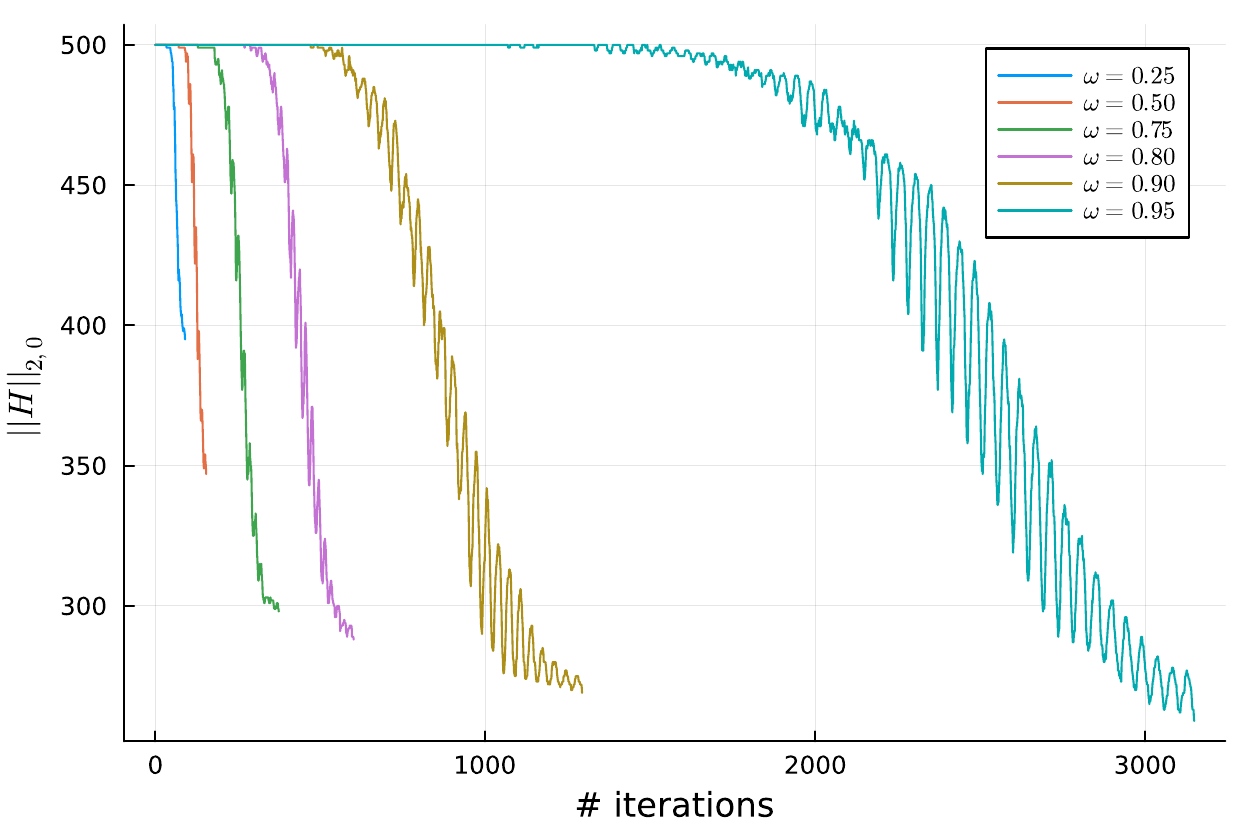} }}%
    ~
    \subfloat[]{{\includegraphics[width=0.48\textwidth]{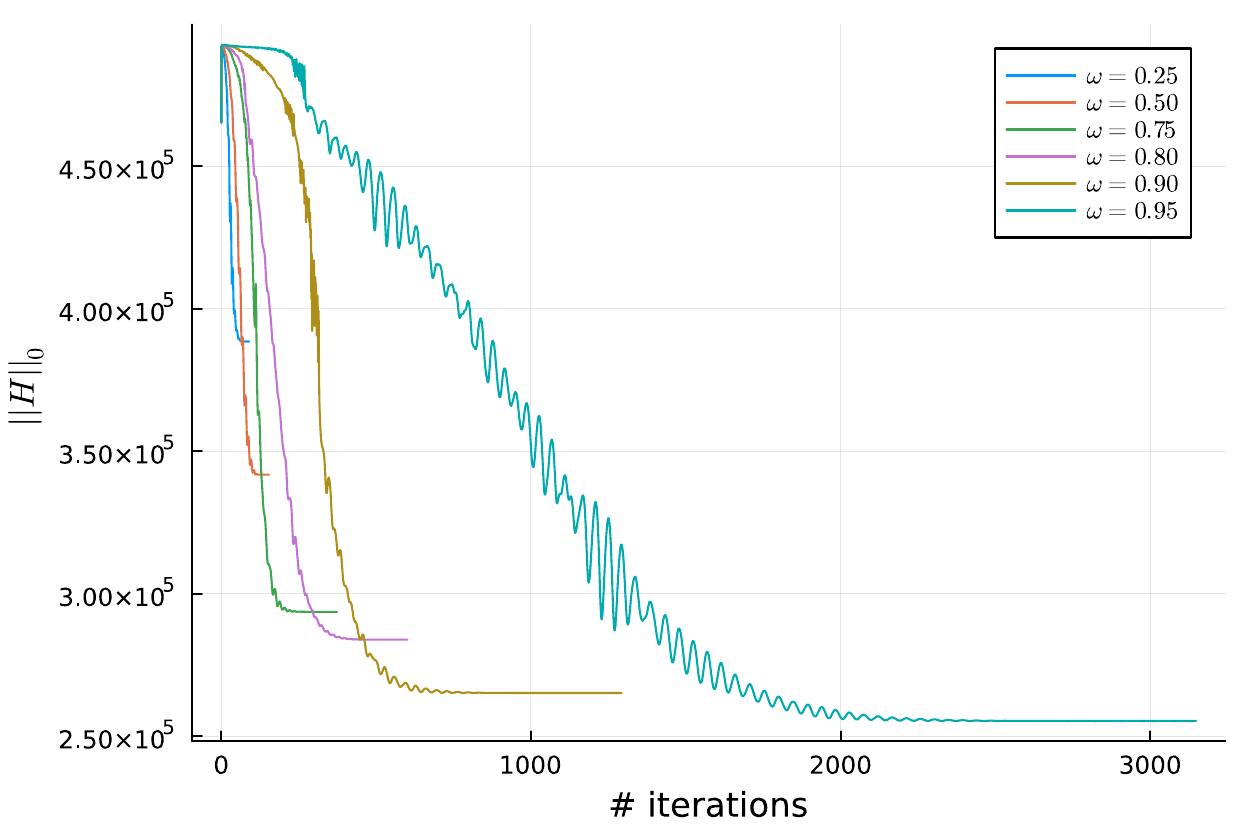} }}%
    \hfill
      \caption{Norms per iteration: ADMM for target 2,0-norm $\omega r + (1-\omega)\|H_{opt}^{2,1}\|_{2,0}$\,; $m,n,r=1000,500,250$ (Instance L1) 
      } 
    \label{fig:compare_bounds3}%
\end{figure}

From the results presented in Table \ref{tab:ADMM2120}, we can  observe  that  \hyperlink{algadmm20norm}{ADMM$_{2,0}$} computes ah-symmetric reflexive generalized inverses with the target  2,0-norms for all instances and  all values of $\omega$ in less than 400 seconds. Increasing $\omega$ leads to slower convergence (in Fig. \ref{fig:compare_bounds3}, we  illustrate this result, where we see a better behavior of the algorithm for $\omega\leq 0.80$).  
Despite disregarding the 1- and 2,1-norms in the formulation of problem \eqref{prob:20converteda}, for which  
\hyperlink{algadmm20norm}{ADMM$_{2,0}$}  is applied, the trade-off between the 2,0-norm and the 1- and 2,1-norms when $\omega$ increases is very clear  from the results in Table \ref{tab:ADMM2120}. The 1- and 2,1-norms of the solutions are relatively small, probably because we start \hyperlink{algadmm20norm}{ADMM$_{2,0}$} with a 
 Frobenius-norm minimizing solution  (see discussion in   \S\ref{subsec:admm1}). This initialization was, in fact, very important in this case. When using different initializations  we observed solutions with  large 1- and 2,1- norms.
For $\omega=0.80$, for example, we see that  
\hyperlink{algadmm20norm}{ADMM$_{2,0}$} is an excellent approach to construct fast structured ah-symmetric reflexive generalized inverses with small 2,1-norms for our test instances, when compared to the 2,1-norm of the local-search procedures, even LS$_{2,1}$\,, that can actually obtain solutions with significant smaller 2,1-norm than LS, but at a high computational cost. We note that \hyperlink{algadmm20norm}{ADMM$_{2,0}$}\,, for $\omega=0.80$, scales even better than the fastest local-search LS, converging in 61 seconds for the largest instance while  LS takes 192.4 seconds to converge. The number of nonzero rows in the solutions obtained with $\omega=0.80$ is always small when compared to $H_{opt}^{2,1}$\,, and the 1- and 2,1-norms are always  small when compared to the local-search solution. Observing now the results for \hyperlink{algadmm2120norm}{ADMM$_{2,1/0}$}\,, we see that it is effective in  obtaining solutions with  smaller  2,1-norm than \hyperlink{algadmm20norm}{ADMM$_{2,0}$}\,, but again, at a much greater computational cost.  An important observation from our experiments is that, with our parameter settings,  the ADMM algorithms proposed for the two nonconvex problems addressed in \S\ref{sec:admm20}, converge  for all tested instances.

\begin{figure}[!ht]
    \centering   \includegraphics[width=0.88\textwidth]{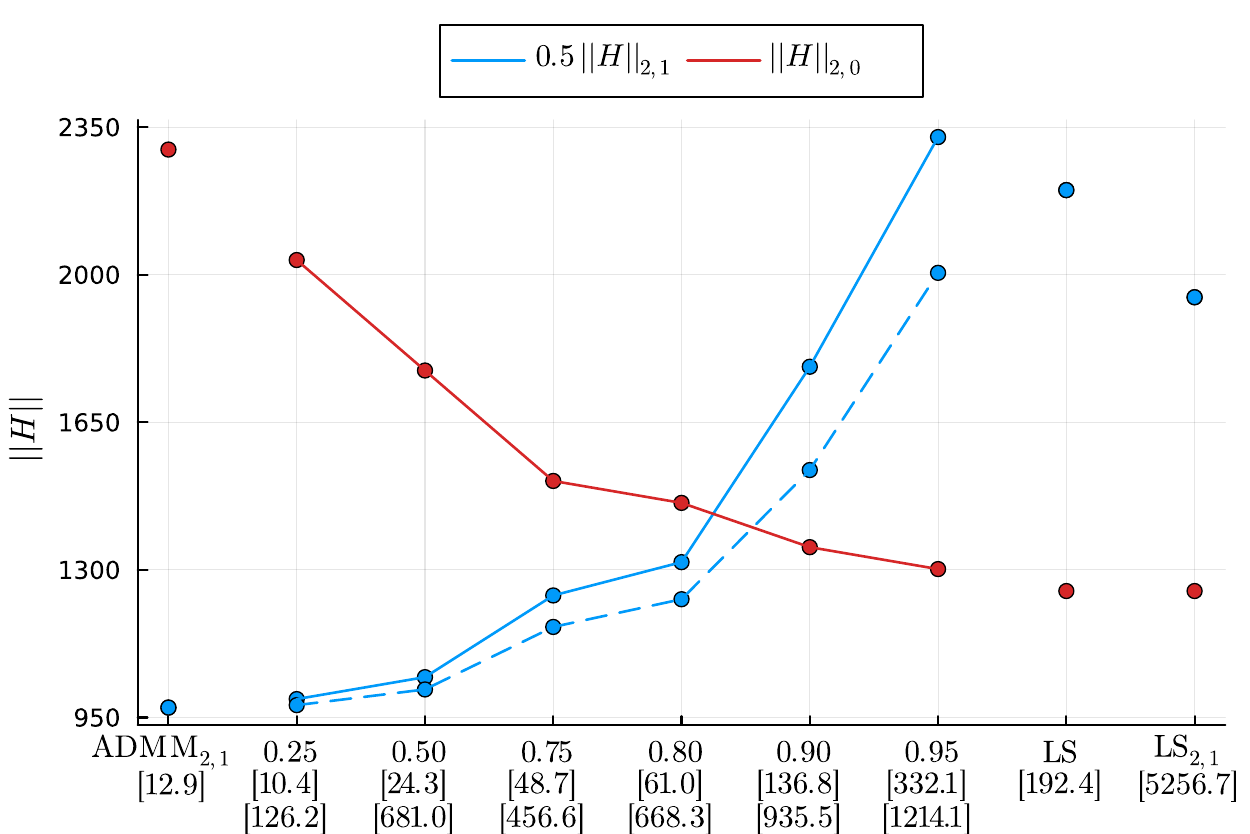}
  \caption{Comparison of our algorithms aimed at row-sparsity on a large instance}\label{fig:20vs2120}
\end{figure}

Fig. \ref{fig:20vs2120} is a visual representation of the
``L5 block'' of Table \ref{tab:ADMM2120}.
We are comparing \hyperlink{algadmm20norm}{ADMM$_{2,0}$}  
and \hyperlink{algadmm2120norm}{ADMM$_{2,1/0}$}\,,  varying $\omega\in\{0.25,0.50,0.75,0.80,0.90,0.95\}$, to obtain a row-sparse solution, for the largest instance that we have considered. Additionally, we compare \hyperlink{algadmm21norm}{ADMM$_{2,1}$} (at the left of the figure), and the local searches  LS and LS$_{2,1}$ (at the right of the figure).
In blue, we show the  2,1-norms of the solutions, scaled by $1/2$ for convenience. The solid blue line corresponds to \hyperlink{algadmm20norm}{ADMM$_{2,0}$} and the dashed blue line corresponds to \hyperlink{algadmm2021norm}{ADMM$_{2,1/0}$}\,. In red, we show the  2,0-norms of the solutions; they are \emph{the same} for  \hyperlink{algadmm20norm}{ADMM$_{2,0}$}  and \hyperlink{algadmm2120norm}{ADMM$_{2,1/0}$}\,. 
Under the horizontal axis, the first row indicates the algorithm, with the numeric values corresponding to $\omega$ for  \hyperlink{algadmm20norm}{ADMM$_{2,0}$}  and \hyperlink{algadmm2120norm}{ADMM$_{2,1/0}$}\,. 
The second row indicates the running times for 
\hyperlink{algadmm21norm}{ADMM$_{2,1}$}\,, \hyperlink{algadmm20norm}{ADMM$_{2,0}$}
for each $\omega$, LS and LS$_{2,1}$\,.
The third row indicates the running times for 
\hyperlink{algadmm2120norm}{ADMM$_{2,1/0}$}\,. 

There is no overall winner, and the best choice  
depends on considering the 2,0-norm, the 2,1-norm, and the elapsed time.
For modest to moderate values of $\omega$, say $\omega\in\{0.25,0.50,0.75,0.80\}$, we can see
a reasonable trade-off in 2,0-norm vs.\! 2,1-norm,
and with a modest but increasing elapsed time.
For the larger values of $\omega$,
we pay a large penalty in the  2,1-norm
for further decrease in the 2,0-norm, and the elapsed time grows as well. Still, some user might prefer the 
solution of say \hyperlink{algadmm2120norm}{ADMM$_{2,1/0}$} at $\omega=0.95$ to the local search LS, giving a 
lesser 2,1-norm and only slightly greater 2,0-norm,
albeit with a much greater elapsed time.
If one is willing to suffer a very long elapsed time,
the local search LS$_{2,1}$ dominates
both the 2,1-norm and  2,0-norm of 
\hyperlink{algadmm2120norm}{ADMM$_{2,1/0}$} at $\omega=0.95$.

\section{Conclusions}\label{sec:Next}

The results from the  numerical experiments performed in this work show that the four ADMM algorithms developed  are  effective in obtaining sparse and row-sparse ah-symmetric reflexive generalized inverses. It is important to note that only two of the four ADMM algorithms address convex optimization problems, and  there is no guarantee of convergence for the other two. Nevertheless, all four of the  algorithms obtain good results and always converge.  

If ensuring convergence to the global optimum is an important point, the algorithms \hyperlink{algadmm1norm}{ADMM$_{1}$} (aimed at inducing sparsity) and \hyperlink{algadmm21norm}{ADMM$_{2,1}$}  (aimed at inducing row-sparsity) are well suited, as they address  convex optimization problems. Both algorithms converge  faster than general-purpose optimization solvers. 
 Comparing both algorithms, the convergence of \hyperlink{algadmm21norm}{ADMM$_{2,1}$} is much faster than \hyperlink{algadmm1norm}{ADMM$_{1}$}\,, and it converges in many fewer iterations. \hyperlink{algadmm21norm}{ADMM$_{2,1}$} is robust and quickly converges to high-precision solutions. As expected, the minimization of the 2,1-norm is more effective than the minimization of the 1-norm if we are aiming for low 2,0-norm (i.e., row-sparse solutions).
 But if we seek low 0-norm (i.e., sparse solutions),
 we naturally should instead minimize the 1-norm.
  ADMM$_{1}^\epsilon$ obtains the  sparsest solutions, but we can expect a high computational cost to have convergence to such solutions. If we are aiming for sparse solutions computed quickly,  then \hyperlink{algadmm21norm}{ADMM$_{2,1}$} and ADMM$_{2,1}^\epsilon$ are preferred to \hyperlink{algadmm1norm}{ADMM$_{1}$}\,.

If we do not require guaranteed convergence to a global optimum, the ADMM algorithms developed for the nonconvex problems addressed, where we limit the number of nonzero rows in the solution, are very effective in obtaining  row-sparse solutions. Moreover, interesting solutions can be obtained with these algorithms by varying the number of nonzero rows allowed. 
Decreasing this number from the 2,0-norm of the 2,1-norm minimizing ah-symmetric reflexive generalized inverse of a given matrix $A$, until it approaches the minimum number of rows (given by $\rank(A)$),  we see an increase in the 1- and 2,1-norms of the solutions. 
With an appropriate number of nonzero rows allowed, \hyperlink{algadmm20norm}{ADMM$_{2,0}$} becomes an excellent approach to quickly construct structured ah-symmetric reflexive generalized inverses with small 2,1-norms. \hyperlink{algadmm20norm}{ADMM$_{2,0}$}   scales even better than our fastest local-search. 
\hyperlink{algadmm2120norm}{ADMM$_{2,1/0}$}\,,  is effective in  obtaining solutions with  smaller  2,1-norm than \hyperlink{algadmm20norm}{ADMM$_{2,0}$}\,, but  at a much greater computational cost.  If one insists on a 
solution with minimum 2,0-norm, then the local searches 
are quite appropriate. In summary, there is no overall winner, and all of our algorithms have their use, as 
one trades off 2,0-norm, 2,1-norm and elapsed time. 

We note that the ADMM algorithms for the 1- and 2,1-norms presented in \S\ref{subsec:admm1}
 and \S\ref{subsec:ADMM21} have their convergence guaranteed  to an optimal solution by general convergence results for ADMM algorithms that can be found, for example, in \cite[Section 3.2]{boyd2011distributed}. The nonconvexity of the feasibility problems addressed in \S\ref{sec:admm20} 
 precludes
 the direct use of these results. 
 There are convergence results for ADMM applied to nonconvex problems, 
  but the ones that we are aware of (for example, \cite{wotao}) do not apply to our situation. Nevertheless, because we were able to get practical convergence
  for our nonconvex ADMMs, exploring a theoretical reason for this
  looks to be a promising direction for future research.

ADMM is just one proximal-gradient algorithm, and there are other such algorithms, and relatives as well. We see our work as opening the door for examining whether any other algorithms of this type may further improve the state-of-the-art for the 
 computation of row-sparse ah-symmetric reflexive generalized inverses.
 


\section*{Acknowledgment}
The authors gratefully acknowledge Laura Balzano and Ahmad Mousavi
for suggesting to us an ADMM approach to structured-sparse generalized inverses, 
but aimed
at an objective minimizing
a balance between the nuclear norm and 2,1-norm of $H$, subject to
\ref{property1} and 
\ref{property3}. By exploiting ideas in \cite{PFLX_ORL},
we  instead worked with \ref{prob:barmin21norm}\,.


\bibliographystyle{plain+eid}


\bibliography{ADMM_LS}


\appendix

\section{Ranks one and two}
\label{app:oneandtwo}
Following the ideas of
\cite{XFLrank12}, we consider the construction of a 2,1-norm
minimizing ah-symmetric reflexive generalized inverse   based  in Theorem \ref{thm:ahconstruction}, for 
the cases where $\rank(A)\in\{1,2\}$. 


\begin{theo}
    Let $A$ be an arbitrary $m \times n$, rank-1 matrix. Choose any row $s \in \{1,\dots,m\}$ (except a zero row), then define $t := \argmax_j\{|A_{sj}|:j = 1,\dots,n\}$. Let $\hat{a}$ be column $t$ of A. Then the $n\times m$ matrix $H$
constructed by Theorem \ref{thm:ahconstruction} over $\hat{a}$, 
is an ah-symmetric reflexive generalized
inverse of A with minimum 2,1-norm.
\end{theo}
\begin{proof}
    First, note that 
    $
    \|\hat{a}^\dagger\|_2 =  \|(\hat{a}^\top \hat{a})^{-1}\hat{a}^\top\|_2 =  \textstyle\left\|\frac{1}{\| \hat{a}\|^2_2}\hat{a}^\top\right\|_2 = \frac{1}{\|\hat{a}\|_2},
    $
and construct  $W$ from Lemma \ref{lem:AWE} with $E = \|\hat{a}^\dagger\|_2$. As we have a rank-1 matrix, $W$ is an $m\times n$ matrix with all elements equal to zero except  $W_{st} = \|\hat{a}^\dagger\|_2/A_{st}$ \,.
Then, $A^\top W$ is a matrix with all elements equal to zero, except for column $t$ which is
$(A^\top W)_{it} =  \textstyle \frac{A_{si}\|\hat{a}^\dagger\|_2}{A_{st}} = \frac{A_{si}}{A_{st} \|\hat{a}\|_2}\,,$
then
$(A^\top W A^\top)_{ij} =  \textstyle \frac{A_{si}A_{jt}}{A_{st}\|\hat{a}\|_2}\,,$
and so 
$\|(A^\top W A^\top)_{i\cdot}\|_2 = \textstyle \frac{|A_{si}|}{|A_{st}|\,\|\hat{a}\|_2}\|\hat{a}\|_2 =  \frac{|A_{si}|}{|A_{st}|}\,,$
for all $i = 1,\dots,n$, $j = 1,\dots,m$.

For $i = t$ we have an active dual constraint, and for $i \neq t$ we have $\|(A^\top W A^\top)_{i\cdot}\|_2 \leq 1$. Then by weak duality, we
have that the constructed $H$ is optimal for \ref{prob:C}.
\end{proof}




  Generally, when $\rank(A) = 2$, we cannot construct a 2,1-norm minimizing ah-symmetric reflexive  generalized inverse based on the column-block construction. Even under the condition that $A$ is totally unimodular and $m = r$, we have the  example:
  $
  A:=\begin{bmatrix}
      1 & 1 & 0\\
      0 & 1 & 1
  \end{bmatrix}.
  $
  We have an ah-symmetric reflexive generalized inverse with minimum 2,1-norm $\frac{\sqrt{2}}{2}\left(1 + \sqrt{3}\right)$,
  $
  H := \textstyle\frac{1}{6}\begin{bmatrix}
 3 + \sqrt{3}~   & \sqrt{3} - 3 \\
  3 - \sqrt{3}~  & 3 - \sqrt{3}\\
  \sqrt{3} - 3~  & 3 + \sqrt{3}
  \end{bmatrix}.
  $
  However, the three ones
  based on our
column block construction have 
2,1-norm $1 +\sqrt{2}$, $2$, $1 + \sqrt{2}$, respectively.

Next, we demonstrate that under an efficiently-checkable technical condition, when $\rank(A) = 2$, construction of a 2,1-norm minimizing ah-symmetric reflexive 
 generalized inverse can be based on the column block construction.
Let $T$ be an ordered subset of $r$ elements from $\{1,\dots,n\}$ and $\hat{A}:=A[\cdot,T]$ be the $m\times r$ submatrix of an $m\times n$ matrix $A$ formed by columns $T$, and $\rank(\hat{A})=r$.
Let $S$ be an ordered subset of $r$ elements from $\{1,\dots,m\}$, such that 
$\tilde{A}:= A[S,T]$ is a  nonsingular matrix. 
\begin{lemm}\label{lem:2normpinvAhat}
Let $k \!\in \!\{1,2\}$ and $\sigma(k)\!:= \!\{1,2\}\!\setminus \!\{k\}$. For rank $2$, 
    $
    \textstyle \|\hat{A}^\dagger_{k\cdot}\|_2^2 = \textstyle \frac{\|\hat{A}_{\cdot\sigma(k)}\|_2^2}{\det(\hat{A}^\top \hat{A})}\,.
    $
\end{lemm}
\begin{proof}
We have $\hat{A}^\dagger_{i\cdot} = (\hat{A}^\top \hat{A})_{i\cdot}^{-1}\hat{A}$, then
$
    \textstyle \|\hat{A}^\dagger_{k\cdot}\|_2^2 = (\hat{A}^\top\hat{A})^{-1}_{k\cdot} \hat{A}^\top \hat{A}  ((\hat{A}^\top\hat{A})^{-1}_{k\cdot})^\top
    = \mathbf{e}_k^\top((\hat{A}^\top\hat{A})^{-1}_{k\cdot})^\top
     = (\hat{A}^\top\hat{A})^{-1}_{k\cdot}\mathbf{e}_k
    = (\hat{A}^\top\hat{A})^{-1}_{kk}\,=  \frac{\|\hat{A}_{\cdot \sigma(k)}\|_2^2}{\det(\hat{A}^\top \hat{A})}\,.
$  
\end{proof}

\begin{theo}
    Let $A$ be an arbitrary $m \!\times \!n$, rank-2 matrix. For any $j_1, j_2 \in \{1,\dots,n\}$, with $j_1 \!<\! j_2$, let $\hat{A} := [
        \hat{a}_{j_1} ~ \hat{a}_{j_2}]$ be the $m\times 2$ submatrix of $A$ formed
by columns $j_1$ and $j_2$. Suppose that $j_1, j_2$ are chosen to minimize the 2,1-norm
of $\hat{H} := \hat{A}^\dagger$ among all $m \times 2$ rank-2 submatrices of $A$. Every column $\hat{b}$ of
$A$, can be uniquely written in the basis $\hat{a}_{j_1},\hat{a}_{j_2}$, say $\hat{b} = \beta_1 \hat{a}_{j_1} +\beta_2 \hat{a}_{j_2}$\,. Suppose
that for each such column $\hat{b}$ of A we have $|\beta_1| + |\beta_2| \leq 1$. Then the $n \times m$ matrix $H$ constructed by Theorem \ref{thm:ahconstruction} based on $\hat{A}$, is an ah-symmetric reflexive generalized inverse of A with minimum 2,1-norm.

\end{theo}

\begin{proof}
    Let $\alpha_1 :=\|\hat{A}^\dagger_{1\cdot}\|_2/{\|\hat{A}_{\cdot 2}\|_2^2}$ and $\alpha_2 :=\|\hat{A}^\dagger_{2\cdot}\|_2/{\|\hat{A}_{\cdot 1}\|_2^2}$\,,
   and construct $W$ from Lemma \ref{lem:AWE} with \[E := \begin{bmatrix}
    \|\hat{A}^\dagger_{1\cdot}\|_2 & -\alpha_1 \hat{A}^\top_{\cdot 1}\hat{A}_{\cdot 2}\\
    -\alpha_2 \hat{A}^\top_{\cdot 1}\hat{A}_{\cdot 2} & \|\hat{A}^\dagger_{2\cdot}\|_2
\end{bmatrix},\]
    so $\hat{A}^\top W A^\top = E\hat{A}^\top$,
and let $k \in \{1,2\}$ and $\sigma(k) := \{1,2\}\setminus \{k\}$. Then we have
\begin{align*}
   & \|(\hat{A}^\top W A^\top)_{k\cdot}\|_2^2 = \|E_{k\cdot} \hat{A}^\top\|_2^2
    =\textstyle \frac{\|\hat{A}^\dagger_{k\cdot}\|_2^2}{\|\hat{A}_{\cdot \sigma(k)}\|_2^4}\left(\|\hat{A}_{\cdot k}\|_2^2\|\hat{A}_{\cdot \sigma(k)}\|_2^4 -(\hat{A}^\top_{\cdot k}\hat{A}_{\cdot \sigma(k)})^2\|\hat{A}_{\cdot \sigma(k)}\|_2^2\right)\\
    &~ =\textstyle\frac{\|\hat{A}^\dagger_{k\cdot}\|_2^2}{\|\hat{A}_{\cdot \sigma(k)}\|_2^2}(\|\hat{A}_{\cdot k}\|_2^2\|\hat{A}_{\cdot \sigma(k)}\|_2^2 -(\hat{A}^\top_{\cdot k}\hat{A}_{\cdot \sigma(k)})^2)\\
    &~=\textstyle\frac{\|\hat{A}^\dagger_{k\cdot}\|_2^2}{\|\hat{A}_{\cdot \sigma(k)}\|_2^2}((\hat{A}^\top \hat{A})_{kk}(\hat{A}^\top \hat{A})_{\sigma(k)\sigma(k)} \!-\!(\hat{A}^\top \hat{A})_{k\sigma(k)}(\hat{A}^\top \hat{A})_{\sigma(k)k})
    =\textstyle \frac{\|\hat{A}^\dagger_{k\cdot}\|_2^2}{\|\hat{A}_{\cdot \sigma(k)}\|_2^2}\det(\hat{A}^\top \hat{A})\!=\!1,
\end{align*}
where the last equation comes from Lemma \ref{lem:2normpinvAhat}. As $\hat{b} = \beta_1 \hat{a}_{j_1} +\beta_2 \hat{a}_{j_2}$\,, then 
\begin{align*}
   \textstyle \|\hat{b}^\top W A^\top\|_2 &= \|\beta_1 \hat{a}_{j_1}^\top W A^\top +\beta_2 \hat{a}_{j_2}^\top W A^\top\|_2 
   ~\leq~ \|\beta_1 \hat{a}_{j_1}^\top W A^\top\|_2 +\|\beta_2 \hat{a}_{j_2}^\top W A^\top\|_2\\
   &= \textstyle |\beta_1|\cdot\| \hat{a}_{j_1}^\top W A^\top\|_2 +|\beta_2|\cdot\| \hat{a}_{j_2}^\top W A^\top\|_2 
   ~=~ |\beta_1| + |\beta_2|\,,
\end{align*}
where the inequality comes from the triangle inequality. Then  $|\beta_1| + |\beta_2| \leq 1 \Rightarrow \|\hat{b}WA^\top\|_2 \leq 1$, so by weak duality, we
establish that the constructed $H$ is optimal.
\end{proof}

\end{document}